\documentclass[12pt,a4paper]{amsart}
\makeatletter
\renewcommand\normalsize{%
    \@setfontsize\normalsize{11.7}{14pt plus .3pt minus .3pt}%
    \abovedisplayskip 10\p@ \@plus4\p@ \@minus4\p@
    \abovedisplayshortskip 6\p@ \@plus2\p@
    \belowdisplayshortskip 6\p@ \@plus2\p@
    \belowdisplayskip \abovedisplayskip}
\renewcommand\small{%
    \@setfontsize\small{9.5}{12\p@ plus .2\p@ minus .2\p@}%
    \abovedisplayskip 8.5\p@ \@plus4\p@ \@minus1\p@
    \belowdisplayskip \abovedisplayskip
    \abovedisplayshortskip \abovedisplayskip
    \belowdisplayshortskip \abovedisplayskip}
\renewcommand\footnotesize{%
    \@setfontsize\footnotesize{8.5}{9.25\p@ plus .1pt minus .1pt}
    \abovedisplayskip 6\p@ \@plus4\p@ \@minus1\p@
    \belowdisplayskip \abovedisplayskip
    \abovedisplayshortskip \abovedisplayskip
    \belowdisplayshortskip \abovedisplayskip}
\ifdefined\pdfpagewidth
\else
\fi
\calclayout
\makeatother

\usepackage{amsmath,amsthm,amssymb,mathtools,amsfonts,mathrsfs}
\usepackage{xspace,xcolor}
\usepackage[breaklinks,colorlinks,citecolor=teal,linkcolor=teal,urlcolor=teal,pagebackref,hyperindex]{hyperref}
\usepackage[alphabetic]{amsrefs}
\usepackage{tikz-cd}
\usepackage[all]{xy}
\usepackage{bbm}
\usepackage{MnSymbol}
\usepackage{stmaryrd}

\usepackage[bbgreekl]{mathbbol}

\DeclareSymbolFontAlphabet{\mathbb}{AMSb}
\DeclareSymbolFontAlphabet{\mathbbl}{bbold}

\newtheorem{theorem}{Theorem}[section]

\newtheorem{proposition}[theorem]{Proposition}
\newtheorem{corollary}[theorem]{Corollary}

\newtheorem{thm}[theorem]{Theorem}
\newtheorem{lemma}[theorem]{Lemma}

\theoremstyle{definition} 
\newtheorem{defn}[theorem]{Definition}
\newtheorem{definition}[theorem]{Definition}

\newtheorem{remark}[theorem]{Remark}

\newcommand{\lie}[1]{\ensuremath{\mathfrak{#1}}}
\newcommand{\g}{{\lie{g}}}

\newcommand{\R}{\bb{R}}

\newcommand{\qu}{/\kern-.7ex/}
\newcommand{\lqu}{\backslash \kern-.7ex \backslash}
\newcommand{\on}{\operatorname} 
\newcommand{\Aut}{\on{Aut}}
 
\newcommand{\Hom}{\on{Hom}}

\title{Orbifold theta functions and mid-age invariants}

\author{Fenglong You}
\address{Department of Mathematics \\ ETH Z\"urich, \\Rämistrasse 101, \\8092 Zürich, \\Switzerland}
\email{fenglong.you@math.ethz.ch}

\thanks{}


\begin{document}
\date{\today}

\begin{abstract} 
We use the orbifold approach to study theta functions in intrinsic mirror symmetry. We introduce a new type of orbifold invariants for snc pairs, called mid-age invariants, and use these invariants to define orbifold invariants associated with the broken line type. Then, we define the orbifold theta functions as generating functions of orbifold invariants with mid-ages. We show that these orbifold theta functions are well-defined and satisfy the multiplication rule. 
\end{abstract}

\maketitle 

\tableofcontents

\section{Introduction}

\subsection{Motivation}

In \cite{You22}, we defined theta functions for smooth log Calabi--Yau pairs in terms of generating functions of two-point relative invariants:
\begin{align}\label{intro-theta-func-def-1}
\vartheta_p=x^{-p}+\sum_{n=1}^{\infty}nN_{p,n}t^{n+p}x^n,
\end{align}
where 
\begin{align}\label{rel-inv-2-pt}
N_{p,n}=\sum_{\beta} \langle [1]_p,[\on{pt}]_n\rangle_{0,2,\beta}^{(X,D)}
\end{align}
for $p\geq 1$. Note that the notation here is slightly different from the notation in \cite{You22} as we switch the order of the first and the second markings.
The goal of this paper is to generalize this to simple normal crossings (snc) pairs. 

The ring of theta functions in the Gross—Siebert program serves as the affine or projective coordinate ring of the mirror construction. Theta functions can either be defined through tropical counting using broken lines or by punctured logarithmic (log) Gromov—Witten invariants associated with the broken line types. The log theta functions are defined locally in \cite{GS21} in terms of two-point punctured log Gromov--Witten invariants of the broken line types and the authors proved that these locally defined theta functions can be patched to global functions.
According to intrinsic mirror symmetry \cite{GS19}, theta functions satisfy the multiplication rule
\[
\vartheta_{\vec p}\star \vartheta_{\vec q}=\sum_{\vec r\in B(\mathbb Z)} N_{\vec p, \vec q, -\vec r} \vartheta_{\vec r}
\]
and the structure constants
\[
N_{\vec p, \vec q, -\vec r}=\sum_{\beta}N_{\vec p, \vec q,-\vec r}^{\beta}t^\beta=\sum_{\beta}\langle [1]_{\vec p},[1]_{\vec q},[\on{pt}]_{-\vec r}\rangle_{0,3,\beta}^{(X,D)}t^\beta
\] 
are punctured log invariants of \cite{ACGS}.

We study it through the orbifold approach \cite{TY20c} to intrinsic mirror symmetry, where we replace punctured invariants by orbifold invariants of multi-root stacks. Theta functions for snc pairs are much more complicated than theta functions (\ref{intro-theta-func-def-1}) for smooth pairs. One may try to replace the invariants in (\ref{rel-inv-2-pt}) by certain two-point invariants of the snc pair. One can see it directly by the virtual dimension constraint that the second marked point can not be of positive contact order to all $D_i$. It needs to be a punctured marked point, that is, the second marked point needs to carry negative (or zero) contact orders to some $D_i$. Naively, one may try to replace $N_{p,n}$ in terms of these two-point invariants with negative contact orders. For example, if $D=D_1+D_2$, one can try to consider the invariants like
\begin{align*}
N_{(p_1,p_2),(n_1,-n_2)}=\sum_{\beta} \langle [1]_{(p_1,p_2)},[\on{pt}]_{(n_1,-n_2)}\rangle_{0,2,\beta}^{(X,D)}t^\beta.
\end{align*}
However, this naive idea in general does not give the theta functions that satisfy the multiplication rule.  

As motivated by \cite{GS21}, theta functions should be defined in terms of invariants associated with the moduli spaces of (log/orbifold) stable maps of a tropical type called the broken line type.  As far as we know, there is no known connection between orbifold Gromov--Witten theory and tropical geometry. These orbifold invariants have not been defined. More generally, the moduli space of stable log maps can be parametrized by types of tropical maps and the discrete data also includes the information of which strata of $D$ that certain irreducible components of the curves are mapped to. However, in orbifold Gromov—Witten theory of multi-root stacks \cite{TY20c}, the discrete data does not include this information. Therefore, it requires a new type of orbifold invariants, which is a tropical refinement of the orbifold invariants in \cite{TY20c}, to capture this information.

Our motivation for studying theta functions is not restricted to the mirror construction in the Gross--Siebert program \cite{GS21}. In the Gross--Siebert program, one mainly focuses on the case when the snc divisor $D$ contains a zero dimensional stratum which is necessary for constructing mirrors. This is an assumption in \cite{GS21}. But we do not make this assumption here and will consider snc pairs in general. For example, a basic case is when the divisor $D$ has two irreducible components: $D=D_1+D_2$. Indeed, there are important motivations for considering general snc pairs that do not necessarily contain a zero dimensional strata. 

Recall that, to construct the intrinsic mirror of a log Calabi--Yau variety, one needs to consider a maximal unipotent monodromy (MUM) degeneration $\pi:\mathscr X\rightarrow S$ and construct the mirror $\mathcal X^\vee_{\on{mum}}$ by taking proj construction to the degree zero part of the relative quantum cohomology of the pair $(\mathscr X_{\on{mum}},\mathscr D_{\on{mum}})$, where $\mathscr D_{\on{mum}}\supset \pi^{-1}(0)$ is a divisor that has a zero dimensional stratum. On the other hand, one can also consider a degeneration $\pi: \mathscr X_{\on{int}}\rightarrow S$ that is not a MUM degeneration. 
The degree zero part of the relative quantum cohomology of the pair $(\mathscr X_{\on{int}}, \mathscr D_{\on{int}})$, where $\mathscr D_{\on{int}}\supset \pi_{\on{int}}^{-1}(0)$ is still well-defined \cite{ACGS}, \cite{GS19}, \cite{TY20c}. One can also construct a space $\mathcal X^\vee_{\on{int}}$ by taking proj construction to it. However, this space is usually of lower dimension, hence, clearly is not the mirror space. A natural question is how does the space $\mathcal X^\vee_{\on{int}}$ relate to the mirror $\mathcal X^\vee_{\on{mum}}$ constructed through a MUM degeneration. 

In the spirit of the Doran--Harder--Thompson conjecture \cite{DHT}, we can consider an intermediate degeneration, that is, a degeneration that is not MUM but can be further degenerated to a MUM degeneration. An example of it, considered in \cite{DHT}, is a Tyurin degeneration of a Calabi--Yau variety X such that the Tyurin degeneration occurs along a locus in the moduli of $X$ that contains a MUM point. 
It was conjectured in \cite{DHT} (see also \cite{DKY21}*{Conjecture 1.6} for a generalization) that there should be a fibration structure on the mirror. For the Tyurin degeneration considered in \cite{DHT}, the base of the fibration is $\mathbb P^1$. For more general degenerations, we expect that the base of the mirror fibration should be $\mathcal X^\vee_{\on{int}}$. In other words, we expect that there is a natural morphism
\[
\pi: \mathcal X^\vee_{\on{mum}} \rightarrow \mathcal X^\vee_{\on{int}},
\]
which is a fibration. Understanding the theta functions of the mirror space $\mathcal X^\vee_{\on{mum}}$ and the theta functions of the base $\mathcal X^\vee_{\on{int}}$ are essential for proving this conjecture. As a future direction, we will use theta functions defined in this paper to study this duality between degenerations and fibrations via intrinsic mirror symmetry.

\subsection{Main results}

The main obstruction of defining theta functions using orbifold invariants is that the orbifold invariants associated with the broken line types have not been defined yet. We need to have the discrete data that can encode the information of how certain irreducible components of the curves are mapped into a specific stratum of $D$. In this paper, we introduce a new type of orbifold invariants called mid-age invariants that exactly capture this information.

When $D$ is smooth, that is, when $D$ has only one irreducible component, orbifold invariants with mid-ages are studied in \cite{You21}. We considered orbifold invariants of root stacks $X_{D,r}$ when the root $r$ is sufficiently large. In \cite{You21}, orbifold invariants of $X_{D,r}$ with a pair of mid-ages mean that there are two markings $p_1$ and $p_2$ such that their ages $k/r$ and $1-k/r$ are mid-ages in the sense that $k,r-k$ are sufficiently larger than $D\cdot \beta$, where $\beta\in H_2(X)$ is the curve class. This assumption ensures that $p_1$ and $p_2$ are mapped into the same rubber space. As proved in \cite{You21}*{Theorem 1.2}, genus zero orbifold invariants with a pair of mid-ages are independent of the value of each age as long as they are sufficiently large. We first generalize this to allow mid-ages $k_a/r, k_b/r$ such that $k_a+k_b=r+c$ where $c\in \mathbb Z$ is not necessarily zero. This allows us to encode the information of the contact order $c$.

\begin{theorem}\label{intro-thm-mid-age-c}[=Theorem \ref{thm-mid-age-c}]
For $r\gg 1$ and a pair of positive integers $\{k_a,k_b\}$ that satisfies 
 \[
k_a+k_b=r+c,
\]
where $c\in \mathbb Z$.
The genus zero cycle class
\[
r^{m_-+1}\tau_*\left(\left[\overline{M}_{0,\vec k, k_a,k_b,\beta}(X_{D,r})\right]^{\on{vir}}\right)\in A_*\left(\overline{M}_{0,m+2,\beta}(X)\times_{X^{m_++m_-+2}}D^{m_++m_-+2}\right)
\] 
of the root stack $X_{D,r}$ is constant in $r$ for sufficiently large $r$. Furthermore, there exists a positive integer $d_0$ such that this cycle class does not depend on the pair $(k_a,k_b)$ as long as $d_0<k_a\leq k_b$.
\end{theorem}
We refer the readers to Theorem \ref{thm-mid-age-c} for the meaning of the notation appearing here.

Following the convention of \cite{TY20c}, the notion of contact orders and ages are used interchangeably. That means a relative marking with positive contact order $k$ corresponds to an orbifold marking with age $k/r$ and a relative marking with negative contact order $k$ corresponds to an orbifold marking with age $1+k/r$. We use contact orders $\mathbf b$ and $-\mathbf b+c$ to denote the mid-ages $\mathbf b/r$ and $1+(-\mathbf b+c)/r$ where $\mathbf b$ stands for a sufficiently large integer.

When $D=D_1+\cdots+D_n$ is snc, where $D_i$ are irreducible components of $D$, we consider orbifold invariants of multi-root stacks. We impose the age conditions such that there is a pair of markings $p_1$ and $p_2$ with mid-ages along $D_i$ for $i\in I\subset \{1,\ldots, n\}$. Then, this gives a way of requiring irreducible components containing markings $p_1$ and $p_2$ to be mapped into the stratum $\cap_{i\in I}D_i$. 


More precisely, we consider the case when the last two markings are mid-ages along some of $D_i$ and have large or small ages (including age $0$) along other $D_i$. By reordering $D_i$, we may assume that the last two markings are of mid-ages along $D_1,\ldots,D_{n^\prime}$, for some $n^\prime\leq n$. We have the following theorem.

\begin{theorem}[=Theorem \ref{thm-mid-age-mul}]\label{intro-thm-mid-age-mul}
 The genus zero cycle class

\[
\left(\prod_{i=1}^nr_i^{m_{i,-}}\prod_{i=1}^{n^\prime}r_i\right)\tau_*\left(\left[\overline{M}_{0,\{\vec k^j\}_{j=1}^m, \vec k_a,\vec k_b,\beta}(X_{D,\vec r})\right]^{\on{vir}}\right)
\]
of the multi-root stack $X_{D,\vec r}$ is constant in $\vec r$ for sufficiently large $\vec r$. Furthermore, there exists $\vec d_0\in \mathbb Q^{n^\prime}$ such that this cycle class does not depend on the pairs $(k_{i,a}, k_{i,b})$ as long as $d_{i,0}<k_{i,a}\leq k_{i,b}$, for $i=1,\ldots, n^\prime$.
\end{theorem}

Theorem \ref{intro-thm-mid-age-mul} holds for orbifold invariants of snc pairs with more than one pair of mid-age markings. We consider it as a way of giving a tropical refinement of orbifold Gromov—Witten invariants. In this paper, we study the orbifold invariants associated with the broken line type. Orbifold invariants associated with stable maps of more general tropical types will be studied elsewhere. The general tropical types are certainly much more complicated.

\begin{remark}
    We use the rooting parameters $\vec r$ and the pair $(\vec k_a,\vec k_b)$ to define mid-age invariants. But mid-age invariants do not depend on these data. Therefore, mid-age invariants are indeed intrinsic to the relative geometry of $(X,D)$. One can also define these invariants without referring to orbifold Gromov--Witten theory. Considering these invariants as orbifold invariants allows us to use the tools from orbifold Gromov--Witten directly.
\end{remark}

\begin{remark}
    In this paper, we mainly focus on the connection of mid-age invariants with intrinsic mirror symmetry. There are actually more reasons (probably more important reasons) to study mid-age invariants. First of all, similar to the mid-age invariants for smooth pairs  studied in \cite{You21}, mid-age invariants for snc pairs play an important role to obtain a loop axiom (an axiom in CohFT) for orbifold invariants of snc pairs. Secondly, higher genus orbifold invariants of multi-root stacks are polynomials in $\vec r$, the constant terms are used to define the formal orbifold invariants of infinite root stacks in \cite{TY20c} (in this paper, we just call them orbifold invariants). Other coefficients of the polynomial are expected to relate to the mid-age invariants. As studied in \cite{You21}, the $r^1$-coefficient of genus one orbifold invariants are genus zero mid-age invariants. Thirdly, mid-age invariants are the main difference between the actual orbifold Gromov--Witten theory of multi-root stacks and the limiting orbifold Gromov--Witten theory in \cite{TY20c}, \cite{FWY}, \cite{FWY19}. Mid-age invariants are the invariants that we do not use when the authors form structural properties of relative/orbifold Gromov--Witten theory in \cite{TY20c}, \cite{FWY}, \cite{FWY19}. Missing these mid-age invariants may be the reason that some other important properties (e.g. Virasoro constraints) have not been proved for the relative/orbifold Gromov--Witten theories in \cite{TY20c}, \cite{FWY}, \cite{FWY19}.
\end{remark}

We now focus on the mid-age invariants associated with the broken line type. We want to encode the information that the first marking is of contact order $\vec p\in B(\mathbb Z)$ and the second marking is of contact order $\vec c$ and the irreducible component containing $p$ mapped into a point in a lowest dimensional stratum $D_{\vec{\mathbf b}}$ of $D$. We replace the second marking with a pair of mid-age markings. The invariants that we will consider are three-point invariants of the following form
\begin{align}\label{inv-mid-age}
\langle [1]_{\vec p},[1]_{\vec{\mathbf b}},[\on{pt}]_{-\vec{\mathbf b}+\vec c}\rangle_{0,3,\beta}^{(X,D)}.
\end{align}
Note that we also impose a point constraint here. If the lowest dimensional stratum $D_{\vec{\mathbf b}}$ is zero dimensional, then the virtual dimension of the moduli space is already zero. Otherwise, we need this point constraint such that the enumerative problem is of dimension zero.

Without loss of generality, we assume that 
\[
    \vec c=(c_1,\ldots, c_l,0,\ldots,0),
    \]
    where $c_i\neq 0$ for $i\in \{1,\ldots, l\}$ and $l \leq n$.
If $D_{\vec c}:=\cap_{i:c_i\neq 0}D_i=\cap_{i=1}^{l}D_i$ is a lowest dimensional stratum of $D$, we take
\[
\vec {\mathbf b}=(\mathbf b_1,\ldots, \mathbf b_{l},0,\ldots, 0),
\]  
where $\mathbf b_i$ are mid-ages along $D_i$, for $i\in\{1,\ldots, l\}$. The choice of mid-ages is unique in this case. 

On the other hand, if $D_{\vec c}$ is not a lowest dimensional stratum, then the choice of a lowest dimensional stratum in $D_{\vec c}$ is not unique. We can consider
    \[
    \vec{\mathbf b}=(\mathbf b_1,\ldots, \mathbf b_s,0,\ldots, 0),
    \]
    where $\mathbf b_i$, for $i\in\{1,\ldots,s\}$, stands for mid-ages along $D_1,\ldots, D_s$. We assume $s > l$ and $\cap_{i=1}^sD_s$ is a lowest dimensional stratum. Without loss of generality, we assume that there is another lowest dimensional stratum  corresponding to a vector of mid-ages of the form:
    \[
    \vec{\mathbf b}^\prime=(\mathbf b_1,\ldots, \mathbf b_{s-1},0,\mathbf b_{s+1},0,\ldots, 0).
    \]
    We prove that the invariant (\ref{inv-mid-age}) is independent of the choice of the lowest dimensional strata:

\begin{theorem}[=Theorem \ref{thm-indep-mid}]
The following identity holds for mid-age invariants
    \begin{align}
    \langle [1]_{\vec p},[1]_{\vec{\mathbf{b}}^\prime},[\on{pt}]_{-\vec{\mathbf{b}}^\prime+\vec c}\rangle_{0,3,\beta}^{(X,D)}=
    \langle [1]_{\vec p},[1]_{\vec{\mathbf{b}}},[\on{pt}]_{-\vec{\mathbf{b}}+\vec c}\rangle_{0,3,\beta}^{(X,D)},
    \end{align}
\end{theorem}

We use these invariants to define the orbifold theta functions as follows.

\begin{definition}[=Definition \ref{def-theta-mid-age}]
For an snc log Calabi--Yau pair $(X,D)$,    the theta functions are defined as follows: Fix $\vec p \in B(\mathbb Z)\setminus \{0\}$. Let $\sigma_{\on{max}}\in \Sigma(X)$ be a maximal cone of $\Sigma(X)$ such that $\mathrm x\in \sigma_{\on{max}}$, then
    \begin{align}\label{intro-iden-def-theta-mid-age}
   \vartheta_{\vec p}(\mathrm x)&:=
   \sum_{\vec k\in {\mathbb Z}^n} \sum_{\beta: D_i\cdot \beta=k_i+p_i}N_{\vec p, \vec{\mathbf b}+\vec k, -\vec{\mathbf{b}}}^\beta t^{\beta} x^{\vec k},
    \end{align}
where $N_{\vec p, \vec{\mathbf b}+\vec k, -\vec{\mathbf{b}}}^\beta$ are orbifold invariants with mid-ages along the divisors $D_{\mathbf b_i}$ such that $D_{\vec {\mathbf b}}=D_{\sigma_{\on{max}}}$ is a lowest dimensional stratum in $D_{\vec k}$; $x^{\vec k}:=x_1^{k_1}\cdots x_n^{k_n}$. If $D_{\vec k}\not\supseteq D_{\sigma_{\on{max}}}$, then $N_{\vec p, \vec{\mathbf b}+\vec k, -\vec{\mathbf{b}}}^\beta=0$.
\end{definition}


We then prove that the orbifold theta functions satisfy the multiplication rule: 
\begin{align}\label{intro-multi-rule-m}
\vartheta_{\vec p}\star \vartheta_{\vec q}=\sum_{\vec r\in B(\mathbb Z)} N_{\vec p, \vec q, -\vec r}\vartheta_{\vec r},
\end{align}
where the structure constants
\[
N_{\vec p, \vec q, -\vec r}:=\sum_{\beta} N_{\vec p, \vec q, -\vec r}^\beta t^\beta=\sum_{\beta} \langle [1]_{\vec p}, [1]_{\vec q}, [\on{pt}]_{-\vec r}\rangle_{0,3,\beta}^{(X,D)}t^\beta
\]
are in terms of orbifold invariants instead of punctured log invariants. 
The multiplication rule (\ref{intro-multi-rule-m}) is a result of the following theorem whose proof is based on the WDVV equation.
\begin{theorem}[=Theorem \ref{thm-wdvv}]
    For $\vec k+\vec l=
\vec s\in \mathbb Z^n$, we have    
\begin{align}
    \sum_{\vec k,\vec l\in \mathbb Z^n}\sum_{\beta_1,\beta_2}N_{\vec p, \vec{\mathbf b}+\vec k, -\vec{\mathbf{b}}}^{\beta_1}N_{\vec q, \vec{\mathbf b}+\vec l, -\vec{\mathbf{b}}}^{\beta_2}=\sum_{\vec r\in B(\mathbb Z),\vec s\in \mathbb Z^n}\sum_{\beta_1,\beta_2} N_{\vec p, \vec q, -\vec r}^{\beta_1}N_{\vec r,\vec{\mathbf b}+\vec s,-\vec{\mathbf b}}^{\beta_2}.
\end{align}
\end{theorem}

\begin{remark}
Besides the fact that we are using different invariants to define structure constants and theta functions, it is also different from \cite{GS21} because \cite{GS21} focuses on the case that the snc divisor $D$ has a zero dimensional stratum. We do not have this assumption here.  \end{remark}

\begin{remark}
    One can try to relate three-point mid-age invariants to certain combinations of two-point invariants without mid-ages via the degeneration formula. However, from the perspective of the multiplicative rule and the WDVV equation (see the proof of Theorem \ref{thm-wdvv}), one naturally gets three-point invariants instead of two-point invariants.
\end{remark}

In Corollary \ref{cor-theta-sm}, we also showed that this definition of theta functions agree with the theta functions in \cite{You22} when $D$ is smooth. 


Following the work of \cite{NR}, \cite{BNR22}, \cite{BNR24}, orbifold Gromov--Witten invariants of root stacks \cite{TY20c} and punctured log Gromov--Witten invariants of \cite{AC}, \cite{GS13}, \cite{ACGS} are different in general. In particular, the structure constants $N_{\vec p, \vec q, -\vec r}$ defined using orbifold invariants and punctured invariants are different in general. 

Following \cite{BNR22} and \cite{BNR24}, these invariants will agree after iterative blow-ups. Recently, \cite{Johnston} showed that orbifold structure constants and log structure constants are equal when the slope sensitivity assumption \cite{BNR22}, \cite{BNR24} is satisfied. 
When the slope sensitivity is satisfied, we also explain in Section \ref{sec:comparison} how these three-point orbifold invariants are related to two-point invariants.

\subsection{Further directions}

\begin{itemize}
    \item Computation of orbifold theta functions. In \cite{You22}, we computed the theta function $\vartheta_1$ for smooth pairs and showed that it is the inverse of the relative mirror map. In an upcoming paper, we will generalize it to snc pairs. It will provide an enumerative meaning of the relative mirror map in \cite{TY20b}. It is expected that the relative mirror map is related to the difference between orbifold and log Gromov—Witten invariants. It was first observed in \cite{You22a}. Recent advances in logarithmic quasimap theory \cite{Shafi23} provides another evidence for this prediction because it is expected that genus zero log quasimap invariants and orbifold quasimap theory should coincide for maximal contact invariants.  On the other hand, orbifold quaismap theory and orbifold Gromov—Witten theory are related by mirror symmetry (quasimap wall-crossing). Therefore, the enumerative meaning of the relative mirror map may provide an explanation of the difference between orbifold and log Gromov—Witten invariants studied in \cite{NR}, \cite{BNR22} and \cite{BNR24}.
    \item Theta functions under degenerations and mirror fibration structures. To under the mirror duality between degenerations and fibrations (the Doran—Harder—Thompson conjecture and its generalization \cite{DKY21}*{Conjecture 1.6}) via the Gross—Siebert program, we will need to understand how theta functions vary under degenerations and how does it give the mirror fibration structures. Our definition of theta functions is explicit, and provides the foundation to this question.  
    \item A tropical refinement of orbifold Gromov--Witten theory. We were told that there is an attempt of the tropicalization of orbifold stable maps by Robert Crumplin. We believe more general mid-age invariants relate to orbifold invariants which are associated with stable maps of certain tropical types. As the broken line type is studied in this paper, we will study the general case elsewhere. It would be interesting to see the relations between the mid-age invariants and Crumplin’s approach.  
    \item Mid-age invariants in higher genus. The approach in this paper works well for genus zero invariants. For higher genus invariants, as studied partially in \cite{You21} for smooth pairs, the mid-age invariants will depend on the mid-age that we choose. Therefore, higher genus invariants require new ideas.
\end{itemize}

\subsection*{Convention} In this paper, when we talk about multi-root stacks $X_{D,\vec r}$, for $\vec r=(r_1,\ldots, r_n)\in (\mathbb Z_{>0})^n$, we always assume all the $r_i$'s are sufficiently large and pairwise coprime. When we say $\vec r$ is sufficiently large, we mean all the $r_i$'s are sufficiently large for $i\in\{1,\ldots, n\}$.


\section*{Acknowledgement}

The author would like to thank Qile Chen, Robert Crumplin, Mark Gross, Samuel Johnston and Navid Nabijou for helpful discussions.

This project has received funding from the European Union’s Horizon 2020 research and innovation programme under the Marie Skłodowska -Curie grant agreement 101025386.

\section{The tropicalization of an snc pair}

\subsection{Notion in tropical geometry}

Let $\mathbf{Cones}$ be the category of rational polyhedral cones. Let $\omega\in \mathbf{Cones}$ be a rational polyhedral cone and $\Lambda_\omega$ be the lattice of integral tangent vectors to $\omega$. Let $\omega_{\mathbb Z}=\omega\cap \Lambda_{\omega}$ be the set of integral points of $\omega$.

Let $G$ be a graph. We write $V(G), E(G)$ and $L(G)$ for the sets of vertices, edges and legs respectively. 
An abstract tropical curve $\omega$ is data $(G,\mathrm g,\ell)$ where
\begin{itemize}
    \item $\mathrm g: V(G)\rightarrow \mathbb N$ is a genus function;
    \item $\ell: E(G)\rightarrow \Hom(\omega_{\mathbb Z},\mathbb N)\setminus \{0\}$ determines edge lengths. 
\end{itemize}
In the genus zero case, we will omit $\mathrm g$. 

Given the data $(G, \mathrm g, \ell)$, we consider the generalized cone complex $\Gamma (G,\ell)$ along with a morphism of cone complexes $\Gamma (G,\ell)\rightarrow \omega$ with fiber over $s\in \on{Int}(\omega)$ being a tropical curve.

Let $X$ be a smooth projective variety and $D\subset X$ a simple normal crossing divisor with $n$ irreducible components:
\[
D=D_1+\cdots+D_n.
\]
For simplicity, we assume that, for any index set $I\subset \{1,\ldots,n\}$, the stratum $D_I:=\cap_{i\in I}D_i$ is connected. Note that $D_I$ can be empty. Let $\on{Div}_D(X)$ be the group of divisors supported on $D$ and $\on{Div}_{D}(X)^*_{\mathbb R}$ be the dual vector space. Let $D_1^*,\ldots, D_n^*$ be the basis dual to $D_1,\dots, D_n$, then the tropicalization $\Sigma(X)$ is a polyhedral cone complex in $\on{Div}_{D}(X)^*_{\mathbb R}$ defined as
\[
\Sigma(X):=\left\{\sum_{i\in I}\mathbb R_{\geq 0}D_i^* | I\subset \{1,\ldots, n\}, D_I\neq \emptyset\right\}.
\]
We denote by $\Sigma_i(X)$ the set of $i$-dimensional cones of $\Sigma(X)$. Let $d_{\sigma,\on{max}}$ be the dimension of the maximal cones of $\Sigma(X)$ containing  the cone $\sigma\in \Sigma(X)$.
Let 
\[
B=\bigcup_{\sigma\in \Sigma(X)}\sigma.
\]
Let $B(\mathbb Z)$ be the set of integer points of $B$.


\begin{definition}
A family of tropical maps to $\Sigma(X)$ over a cone $\omega$ is morphism of cone complexes
\[
h: \Gamma (G,\ell)\rightarrow \Sigma(X).
\]
\end{definition}

The type of a family of tropical maps consists of data $\tau:=(G, \mathbf \sigma, \mathbf u)$, where
\begin{itemize}
    \item \[
\mathbf \sigma=V(G)\cup E(G)\cup L(G)\rightarrow \Sigma(X)
\]
maps $x\in V(G)\cup E(G)\cup L(G)$ to the minimal cone $\tau\in \Sigma(X)$ containing $h(\omega_x)$.
\item For each edge or leg $E\in E(G)\cup L(G)$, we associate a contact order $\mathbf u(E)\in \Lambda_{\sigma(E)}$, where $\Lambda_{\sigma(E)}$ is the image of the tangent vector $(0,1)\in \Lambda_{\omega_E}=\Lambda_\omega\oplus \mathbb Z$ under the map $h$. 
\end{itemize}
\if{
A type $\tau$ is realizable if there exists a family of tropical maps to $\Sigma(X)$ of type $\tau$. A type $\tau$ is balanced if, for each $v\in V(G)$ with $\sigma(v)\in \Sigma(X)$ a codimension zero or one cone, the balancing condition holds at $v$:
\[
\sum u(E_i)=0,
\]
where $E_i$'s are legs and edges adjacent to $v$ and $u(E_i)$ are contact orders.
}\fi
Following \cite{GS21}*{Definition 3.19}, a non-trivial broken line type $\tau=(G, \mathbf \sigma, \mathbf u)$ of tropical map to $\Sigma(X)$ is realizable and balanced such that:
\begin{itemize}
    \item $G$ is a genus zero graph with $L(G)=\{L_{\on{in}}, L_{\on{out}}\}$ with $\sigma(L_{\on{out}})\in \Sigma(X)$ and
    \[
    u_{\tau}:=\mathbf u(L_{\on{out}})\neq 0, \quad p_\tau:=\mathbf u(L_{\on{in}})\in \sigma(L_{\on{in}})\setminus \{0\}.
    \]
    \item Let $h:\Gamma(G,\ell) \rightarrow \Sigma(X)$ be the corresponding universal family of tropical maps, and let $\tau_{\on{out}}\in \Gamma(G,\ell)$ be the cone corresponding to $L_{\on{out}}$. Then $\dim \tau=d_{\sigma, \on{max}}-1$ and $\dim h(\tau_{\on{out}})=d_{\sigma,\on{max}}$.
\end{itemize}

Under the assumption of \cite{GS21}, the virtual dimension of the moduli space is zero, no additional constraint is needed. When $D$ does not contain a zero dimensional stratum, we also impose a point constraint at $L_{\on{out}}$.

A decorated type is data ${\boldsymbol \tau}=(G,\mathbf \sigma, \mathbf u, \mathbf A)$ where $\mathbf A: V(G)\rightarrow H_2(X)$ associates a curve class to each vertex of $G$. The total curve class of $\mathbf A$ is $A=\sum_{v\in V(G)}\mathbf A(v)$.


\section{Preliminary on orbifold Gromov--Witten theory}
We briefly review the formal Gromov--Witten theory of infinite root stacks as a limit of the orbifold Gromov--Witten theory of multi-root stacks. The foundation of this orbifold Gromov--Witten theory can be found in \cite{TY20c}. We refer readers to  \cite{AV}, \cite{AGV02}, \cite{AGV}, \cite{CR01} and \cite{Tseng} for the foundation of orbifold Gromov--Witten theory in general.

\subsection{The orbifold pairing}
Let $\mathcal X$ be a Deligne--Mumford stack, the inertia stack $\mathcal {IX}$ is defined to be the fiber product 
\[
\mathcal {IX}:=\mathcal X \times_{\Delta,\mathcal X\times \mathcal X, \Delta}\mathcal X,
\]
where $\Delta: \mathcal X\times \rightarrow \mathcal X\times \mathcal X$ is the diagonal morphism. The objects of $\mathcal {IX}$ are the following
\[
\on{OB}(\mathcal {IX})=\{(x,g)| x\in \on{OB}(\mathcal X), g\in \Aut_{\mathcal X}(x)\}.
\]
Let
\[
\mathcal {IX}=\coprod_{i\in \mathcal I}\mathcal X_i,
\]
where $\mathcal I$ is an index set. There is a natural involution map 
\[
I:\mathcal {IX}\rightarrow \mathcal {IX}. 
\]
On objects, we have 
\[
I((x,g))=(x,g^{-1}).
\]
Restricting to $\mathcal X_i$, the map $I$ is an isomorphism between $\mathcal X_i$ and another component, denoted by $\mathcal X_{i^{I}}$. The orbifold Poincar\'e pairing is defined as follows. For $\alpha_1 \in \mathcal X_i$ and $\alpha_2\in \mathcal X_{i^I}$,
\[
\int_{\overline{I}\mathcal X}\alpha_1\cup \alpha_2:=\int_{\mathcal X_i}\alpha_1\cup I^*\alpha_2.
\]

\subsection{Orbifold Gromov--Witten theory of multi-root stacks}\label{sec:GW-multi-root}
We consider orbifold Gromov--Witten theory of multi-root stacks. The limit is called the formal Gromov--Witten theory of infinite root stacks in \cite{TY20c}.

Let $X$ be a smooth projective variety and
\[
D=D_1+\ldots+D_n
\]
be an snc divisor with irreducible components $D_1,\ldots, D_n$. For any index set $I\subseteq\{1,\ldots, n\}$, we define 
\[
D_I:=\cap_{i\in I} D_i.
\]
Let 
\[
\vec k=(k_1,\ldots,k_n)\in \mathbb Z^n.
\]
We define
\[
I_{\vec k}:=\{i:k_i\neq 0\}\subseteq \{1,\ldots,n\}.
\]

For $\vec r=(r_1,\ldots,r_n)\in (\mathbb Z_{\geq 0})^n$ and $r_i$'s are pairwise coprime, we consider the multi-root stack
\[
X_{D,\vec r}:=X_{(D_1,r_1),\ldots,(D_n,r_n)}.
\]
 By \cite{TY20c}*{Corollary 16}, genus zero orbifold Gromov-Witten invariants of $X_{D,\vec r}$, after multiplying suitable powers of $r_i$, is independent of $r_i$ for $r_i$ sufficiently large. Following \cite{TY20c}*{Definition 18}, the formal Gromov-Witten invariants of $X_{D,\infty}$ are defined as limits of the corresponding genus zero orbifold Gromov-Witten invariants of $X_{D,\vec r}$.

Let \begin{itemize}
    \item $\gamma_j\in H^*(D_{I_{\vec k^j}})$, for $j\in\{1,2,\ldots,m\}$;
    \item $a_j\in \mathbb Z_{\geq 0}$, for $j\in \{1,2,\ldots,m\}$.
\end{itemize}
The formal genus zero Gromov-Witten invariants of $X_{D,\infty}$ are defined as
\begin{align}\label{inv-X-D}
    \left\langle [\gamma_1]_{\vec k^1}\bar{\psi}^{a_1},\ldots, [\gamma_m]_{\vec k^m}\bar{\psi}^{a_m} \right\rangle_{0,\{\vec k^j\}_{j=1}^m,\beta}^{X_{D,\infty}}:=\left(\prod_{i=1}^n r_i^{k_{i,-}}\right)\left\langle \gamma_1\bar{\psi}^{a_1},\ldots, \gamma_m\bar{\psi}^{a_m} \right\rangle_{0,\{\vec k^j\}_{j=1}^m,\beta}^{X_{D,\vec r}}
\end{align}
for sufficiently large $\vec r$, where
\begin{itemize}
    \item the vectors
\[
\vec k^j=(k_1^j,\ldots,k_n^j)\in (\mathbb Z)^n, \text{ for } j=1,2,\ldots,m,
\]
satisfy the following condition:
\[
\sum_{j=1}^m k_i^j=\int_\beta[D_i], \text{ for } i\in\{1,\ldots,n\}.
\]
We also use the notation
\[
\sum_{j=1}^m\vec k^j=D\cdot \beta.
\]
These vectors are used to record contact orders of the markings with respect to divisors $D_1,\ldots,D_n$. In the language of orbifold Gromov--Witten theory,
the $j$-th marking maps to twisted sector $D_{I_{\vec s^j}}$ of the inertia stack of $X_{D,\vec r}$ with age
\[
\sum_{i: k^j_i>0} \frac{k^j_i}{r_i}+\sum_{i: k^j_i<0} \left(1+\frac{k^j_i}{r_i}\right),
\]
see \cite{TY20c}. 
\item Define
\[
k_{i,-}:=\#\{j: k_i^j<0\}, \text{ for } i=1,2,\ldots, n,
\]
to be the number of markings that have negative contact order with the divisor $D_i$.
\end{itemize}

The genus zero invariant (\ref{inv-X-D}) is zero unless it satisfies the virtual dimension constraint
\begin{align}\label{vir-dim}
    \dim_{\mathbb C}X-3+m+\int_\beta c_1(T_X)-\int_\beta [D]=\sum_{j=1}^m \deg^{0}([\gamma_j]_{\vec k^j})+\sum_{j=1}^m a_j,
\end{align}
where if $\alpha\in \mathfrak H_{\vec k}$ is a cohomology class of real degree $d$, then,
\begin{align}\label{deg-0}
\deg^{0}([\alpha]_{\vec k})=d/2+\#\{i:k_i<0\}.
\end{align}
\begin{remark}
Besides Section \ref{sec:comparison}, we only consider orbifold invariants not log invariants. Therefore, we will sometimes use the notation  $\langle \cdots \rangle^{(X,D)}$ instead of $\langle \cdots\rangle^{X_{D,\infty}}$ when there is no ambiguity. As the contact orders $\{\vec k^j\}$ are specified in the insertions, we will also use the notation $\langle \cdots \rangle_{0,m,\beta}^{(X,D)}$ instead of $\langle \cdots \rangle_{0,\{\vec k^j\}_{j=1}^m,\beta}^{(X,D)}$.
\end{remark}

\section{The degeneration formula for orbifold Gromov--Witten theory}

\subsection{General case}

We recall the degeneration formula for orbifold Gromov--Witten theory from \cite{AF}, then we will specialize to the case of multi-root stacks. Note that, we only consider degenerations where the central fiber contains two irreducible components intersecting along a smooth orbifold divisor. 


\begin{defn}[Definition 4.6, \cite{Li01}]
{\emph {An admissible graph}} $\Gamma$ is a graph without edges plus the following data.
\begin{enumerate}
    \item An ordered collection of legs.
    \item An ordered collection of weighted roots.
    \item A function $\g:V(\Gamma)\rightarrow \mathbb Z_{\geq 0}$.
    \item A function $b:V(\Gamma)\rightarrow H_2(X,\mathbb Z)$.
\end{enumerate}
\end{defn}

\if{
\begin{defn}\label{def:admgraph0}
\emph{A (connected) graph of type $0$} is a weighted graph $\Gamma^0$ consisting of a
single vertex, no edges, and the following.
\begin{enumerate}
\item {$0$-roots},
\item {$\infty$-roots of node type},
\item {$\infty$-roots of marking type},
\item {Legs}.
\end{enumerate}
$0$-roots are weighted by
positive rational numbers, and $\infty$-roots are weighted by negative rational numbers. The vertex is associated with a tuple $(g,\beta)$ where $g\geq 0$
and $\beta\in H_2(D,\mathbb Z)$. 
\end{defn}

A graph $\Gamma^\infty$ of type $\infty$ is an admissible graph such that the roots are distinguished by node type and marking type.

\begin{definition}\label{defn:locgraph}
\emph{An admissible bipartite graph} $\mathfrak G$ is a tuple
$(\mathfrak S_0,\Gamma^\infty,I,E,g,b)$, where 
\begin{enumerate}
\item $\mathfrak S_0=\{\Gamma_i^0\}$ is a set of graphs of type $0$;
  $\Gamma^\infty$ is a (possibly disconnected) graph of type $\infty$.
\item $E$ is a set of edges.
\item $I$ is the set of markings.
\item $g$ and $b$ represent the genus and the degree respectively.
\end{enumerate}
Moreover, the admissible bipartite graph must satisfy some conditions, see, for example, \cite{FWY}*{Definition 4.8}.
\end{definition}
}\fi

Let $\mathcal X$ be a smooth proper Deligne--Mumford stack with a projective coarse moduli space $X$. We consider a degeneration of $\mathcal X$, a one-parameter family with a smooth total space $W$, such that the generic fiber is $\mathcal X$ and the central fiber is
\[
W_0=\mathcal X_1\sqcup_{\mathcal D} \mathcal X_2
\]
with $\mathcal X_1$, $\mathcal X_2$ and $\mathcal D$ are smooth Deligne--Mumford stacks.
The (numerical version of the) degeneration formula is the following
\begin{align}
    \langle \prod_{i=1}^m\tau_{a_i}(\gamma_i)\rangle_{g,m,\beta}^{\mathcal X}
    = \sum \frac{\prod_i \eta_i}{\Aut(\eta)}\langle \prod_{i\in S_1}\tau_{m_i}(\gamma_i)| \eta\rangle_{g_1,\eta,|S_1|,\beta_1}^{\bullet, (\mathcal X_1,\mathcal D)}\langle \eta^\vee | \prod_{i\in S_2}\tau_{m_i}(\gamma_i)\rangle_{g_2,\eta,|S_2|,\beta_2}^{\bullet,(\mathcal X_2,\mathcal D)}. 
\end{align}

Notation on the LHS,
\begin{itemize}
    \item $g$, the genus of the curve, is a non-negative integer and $\beta$ is a curve class on $\mathcal X$. 
    \item $\gamma_1,\ldots, \gamma_m \in H^*_{\on{CR}}(\mathcal X):=H^*(\mathcal{IX})$ are cohomology classes.
    \item $a_1,\ldots, a_m$ are non-negative integers.
\end{itemize}

Notation on the RHS,
\begin{itemize}
    \item \[
    \eta=\{(\eta_1,\delta_1),\ldots, (\eta_l, \delta_l)\}
    \]
    is a orbifold cohomology weighted partition with a partition
    \[
    \sum_{i=1}^l \eta_i=\mathcal D\cdot \beta
    \]
    and cohomology classes $\sigma_1,\ldots, \sigma_l\in H^*_{\on{CR}}(\mathcal D)$. Note that $\eta_i$ are non-negative rational numbers.
    \item $\eta^\vee$ is the dual cohomology weighted partition: 
    \[
    \eta^\vee=\{(\eta_1,\delta_1^\vee),\ldots, (\eta_l, \delta_l^\vee)\},
    \]
    where cohomology classes $\delta_1^\vee,\ldots, \delta_l^\vee\in H^*_{\on{CR}}(\mathcal D)$. are Poincar\'e dual to $\delta_1,\ldots, \delta_l\in H^*_{\on{CR}}(\mathcal D)$:
    \begin{align}\label{orb-pairing}
    \int_{\mathcal {\overline{I}D}}\delta_i\cup \delta_i^\vee=1.
    \end{align}
    \item The sign $\bullet$ means possibly disconnected invariants.
\item $S_1\coprod S_2=\{1,\ldots, m\}$ is a splitting of the markings.
\item The admissible bipartite graph $\Gamma$ obtained by gluing admissible graphs $\Gamma_1$ and $\Gamma_2$ along roots is connected of genus $g$ and total weight $\mathcal D\cdot \beta$:
\begin{align}\label{iden-g}
    g=\sum_v g(v)+|E(\Gamma)|-|V(\Gamma)|+1,
\end{align}
\[
\mathcal D\cdot \beta=\sum_v d_v.
\]

    \item The summation is over  all intermediate cohomology weighted partitions and all possible splitting of $g, \beta$ and $m$.
\end{itemize}



\subsection{The degeneration formula for multi-root stacks}\label{sec:deg}

We then specialize to the degeneration formula for genus zero Gromov--Witten invariants of multi-root stacks when the degeneration is the degeneration to the normal cone of one of the divisors $D_i$. Recall that we consider an $m\times n$ matrix of integers
\[
\{k_{i}^j\}_{1 \leq i \leq n, 1\leq j \leq m}
\]
that satisfies
\[
\sum_{j=1}^m k_{i}^j=D_j\cdot \beta \quad \text{ for } 1\leq i \leq n.
\]
Note that we choose $k_{i}^j$ independently from $\vec r$. We consider the genus zero degree $\beta\in H_2(X)$, $m$-pointed orbifold Gromov--Witten invariant of $\mathcal X:=X_{D,\vec r}$, where the insertions are $\gamma_j\in H^*_{\on{CR}}(X_{D,\vec r})$ and the ages along markings $p_j$ are
\[
\on{age}_{p_i}(N_{\mathcal D_j/\mathcal X})=\left\{\begin{array}{cc}
    \frac{k_{i}^j}{r_i}, &  k_{i}^j\geq 0;\\
    1+\frac{k_{i}^j}{r_i}, &  k_{i}^j <0,
\end{array}\right.
\] 
where $\mathcal D_i:=\sqrt[r_i]{N_{D/X}}$ is a gerbe over $D$ banded by $\mu_{r_i}$.
For $\vec r$ sufficiently large, we also refer to $k_{i}^j$ as the contact order of the marking $p_j$ along $D_i$.

We write down the degeneration formula when $D=D_1+D_2$ has two irreducible components. The general case with $n$ irreducible components works similarly.

We consider the multi-root stack $X_{D,\vec r}:=X_{(D_1,r_1),(D_2,r_2)}$ with $r_1$-th root along $D_1$ and $r_2$-th along $D_2$. By \cite{TY20c}*{Definition 16}, the genus zero invariants of the multi-root stack $X_{D,\vec r}$ are the same as the relative orbifold invariants of $(X_{D_2,r_2},D_{1,(D_{12},r_2)})$ when $r_1$ is sufficiently large, where $D_{12}:=D_1\cap D_2$ and $D_{1,(D_{12},r_2)}$ is the $r_2$-th root stack of $D_1$ along $D_{12}$. We consider the degeneration to the normal cone to $D_1$ and apply the orbifold degeneration formula to study the relative orbifold invariants of $(X_{D_2,r_2},D_{1,(D_{12},r_2)})$ (with negative contact orders). 
Let $\mathcal D_1\subset X_{D,\vec r}$ be the $r_1$-th root gerbe of $D_0:=D_{1,(D_{12},r_2)}$. 
The degeneration formula can be written as follows:
\begin{align}
    \langle \prod_{i=1}^m\tau_{a_i}(\gamma_i)\rangle_{g,m,\beta}^{X_{D,\vec r}}
    = \sum \frac{\prod_i \eta_i}{\Aut(\eta)}\langle \prod_{i\in S_1}\tau_{m_i}(\gamma_i)| \eta\rangle_{g_1,\eta,|S_1|,\beta_1}^{\bullet, (X_{D,\vec r},\mathcal D_1)}\langle \eta^\vee | \prod_{i\in S_2}\tau_{m_i}(\gamma_i)\rangle_{g_2,\eta,|S_2|,\beta_2}^{\bullet,(\mathbb P_{\mathcal D_1}(\mathcal O\oplus N_{-\mathcal D_1}),\mathcal D_\infty)}, 
\end{align}
where $\mathcal D_{0,r_1}$ and $\mathcal D_\infty$ are zero and infinity divisors of the bundle $\mathbb P_{\mathcal D_1}(\mathcal O\oplus N_{-\mathcal D_1})$. Note that $\mathcal D_\infty=D_{1,(D_{12},r_2)}$ is the $r_2$-th root stack of $D_1$ along $D_{12}$ and $\mathcal D_{0,r_1}$ is the $r_1$-th root gerbe of $D_0:=D_{1,(D_{12},r_2)}$. Note that, when applying the degeneration formula to $(X_{D_2,r_2},D_{1,(D_{12},r_2)})$, the cohomology weighted partition $\eta$ also carries orbifold ages along $D_{12}$.

For $\vec r$ sufficiently large, genus zero Gromov--Witten theory of $(\mathbb P_{\mathcal D_1}(\mathcal O\oplus N_{-\mathcal D_1}),\mathcal D_\infty)$ equals to the genus zero orbifold Gromov--Witten theory of $\mathbb P_{D_{1,(D_{12},r_2)}}(\mathcal O\oplus N_{-\mathcal D_1})_{(D_0+D_\infty,(r_1,r_1))}$. To simplify the notation slightly, we simply write this pair as $(\mathbb P_{\mathcal D_1}(\mathcal O\oplus N_{-\mathcal D_1}),\mathcal D_0+\mathcal D_\infty)$. We can write the degeneration formula as
\begin{align}\label{degen-formula-orb}
\langle \prod_{i=1}^m\tau_{a_i}(\gamma_i)\rangle_{g,m,\beta}^{(X,D)}
    = \sum \frac{\prod_i \eta_i}{\Aut(\eta)}\langle \prod_{i\in S_1}\tau_{m_i}(\gamma_i)| \eta\rangle_{g_1,\eta,|S_1|,\beta_1}^{\bullet, (X,D)}\langle \eta^\vee | \prod_{i\in S_2}\tau_{m_i}(\gamma_i)\rangle_{g_2,\eta,|S_2|,\beta_2}^{\bullet,(\mathbb P_{\mathcal D_1}(\mathcal O\oplus N_{-\mathcal D_1}),\mathcal D_0+\mathcal D_\infty)}.
\end{align}
Note that the invariants in (\ref{degen-formula-orb}) are orbifold invariants of root stacks with sufficiently large roots, not punctured logarithmic invariants.
Here we need to explain the contact orders (along $D_1$) and ages (along $D_2$) of $\eta$ and $\eta^\vee$. Since $r_1,r_2$ are sufficiently large, contact orders and ages are used interchangeably via the usual relative/orbifold correspondence \cite{ACW}, \cite{TY18}, \cite{FWY}, \cite{FWY19} and \cite{TY20c}. If the $i$-th marking of $\eta$ has contact order
\[
(\eta_{1}^i,\eta_{2}^i),
\]
then the $i$-th marking of $\eta^\vee$ has contact order
\[
(\eta_{1}^i,-\eta_{2}^i).
\]


Note that the first components (relative divisor of the degeneration) have the same contact order as usual; the second components (orbifold divisor of the degeneration) have the opposite contact order (age) because of the orbifold condition (\ref{orb-pairing}) for the degeneration formula for orbifold Gromov--Witten theory. In general, if $D$ has $n$ irreducible components and we consider the degeneration to the normal cone of $D_1$, we can again apply the orbifold degeneration formula. 
If the $i$-th marking of $\eta$ has contact order
\[
(\eta_{1}^i,\eta_{2}^i,\ldots, \eta_{n}^i),
\]
then the $i$-th marking of $\eta^\vee$ has contact order
\[
(\eta_{1}^i,-\eta_{2}^i,\ldots, -\eta_{n}^i).
\]

\subsection{Special cases}
We would like to specialize to the case of orbifold invariants with only one large age. We first need some results in orbifold Gromov--Witten theory. Let $\mathcal X$ be a smooth Deligne--Mumford stack with projective coarse moduli space $X$. Let $\mathcal D\subset \mathcal X$ be a smooth irreducible divisor and $L:=N_{\mathcal D/\mathcal X}$ be the normal bundle of $\mathcal D$ in $\mathcal X$.  We recall that $\sqrt[r]{L/\mathcal D}$ is a gerbe over $\mathcal D$ banded by $\mu_r$. 

Let 
\[
\vec k=(k_1,\ldots,k_m)\in \mathbb Q^m
\]
such that
\[
\sum_{i=1}^m k_i=c_1(L)\cdot \beta, \text{ for } \beta \in H_2(D).
\]
Let $\overline{M}_{0,\vec k}(\sqrt[r]{L/\mathcal D},\beta)$ be the moduli space of genus zero orbifold stable maps of degree $\beta\in H_2(X)$ to $\sqrt[r]{L/\mathcal D}$ with
\[
\on{age}_{p_i}(L)=\left\{\begin{array}{cc}
 \frac{k_i}{r}    & k_i\geq 0 \\
 \frac{k_i}{r}+1    & k_i<0. 
\end{array}\right.
\]
Let $\overline{M}_{0,\vec k}^\sim(\mathcal D,\beta)$ be the moduli space of relative stable maps to rubber targets over $\mathcal D$ with contact orders $\{k_i\}_{i:k_i>0}$ and $\{-k_i\}_{i:k_i<0}$ at $0$ and infinity divisors respectively; and there are $m_0:=\#\{i:k_i=0\}$ interior markings. We have the following forgetful maps
\[
\tau_1: \overline{M}_{0,\vec k}(\sqrt[r]{L/\mathcal D},\beta)\rightarrow \overline{M}_{0,m}(\mathcal D,\beta)
\]
and
\[
\tau_2: \overline{M}_{0,\vec k}^\sim(\mathcal D,\beta)\rightarrow \overline{M}_{0,m}(\mathcal D,\beta).
\]

We have the following lemma.

\begin{lemma}\label{lemma-gerbe-rubber}
   For $r$ sufficiently large, we have
    \[
    r(\tau_1)_*\left( [\overline{M}_{0,\vec k}(\sqrt[r]{L/\mathcal D},\beta)]^{\on{vir}}\right)=(\tau_2)_*\left([\overline{M}_{0,\vec k}^\sim(\mathcal D,\beta)]^{\on{vir}}\right)=[\overline{M}_{0,m}(\mathcal D,\beta)]^{\on{vir}}.
    \]
\end{lemma}
\begin{proof}
Let $$\mathcal Y:=\mathbb{P}(N_{\mathcal D/\mathcal X}\oplus\mathcal{O}_\mathcal X)\xrightarrow{\pi} \mathcal D.$$
The infinity and zero divisors of $\mathcal Y$ are
$$\mathcal D_\infty:=\mathbb{P}(N_{\mathcal D/\mathcal X}) \subset \mathbb{P}(N_{\mathcal D/\mathcal X}\oplus\mathcal{O}_\mathcal X);$$
$$\mathcal D_0:= \mathbb{P}(\mathcal{O}_{\mathcal X})\subset  \mathbb{P}(N_{\mathcal D/\mathcal X}\oplus\mathcal{O}_{\mathcal X}).$$
We consider the $r$-th root stack of $\mathcal Y$ along $\mathcal D_0$:
\[
\mathcal Y_{\mathcal D_0, r}:=\mathbb{P}(N_{\mathcal D/\mathcal X}\oplus \mathcal{O}_\mathcal X)_{\mathcal D_0, r},
\]
and study the orbifold-relative Gromov-Witten invariants of $(\mathcal Y_{\mathcal D_0, r}, \mathcal D_\infty)$. Let $\overline{M}_{0, \vec k}(\mathcal Y_{\mathcal D_0,r}/\mathcal D_\infty,\beta)$ be the moduli space of genus zero degree $\beta$ orbifold-relative stable maps to $(\mathcal Y_{\mathcal D_0, r}, \mathcal D_\infty)$ where, for $k_i>0$, the $i$-th marking is an orbifold marking that maps to $\mathcal D_0$ with age $\on{age}_{p_i}(L)=k_i/r$; for $k_i<0$, the $i$-th marking is a relative marking that maps to $\mathcal D_\infty$ with contact order $-k_i$; for $k_i=0$, the $i$-th marking is an interior marking (which maps to $\underline{I}\mathcal Y$.)

We would like to compute
\begin{align}\label{push-forward-cycle}
(\tau^\prime)_*\left[\overline{M}_{0, \vec k}(\mathcal Y_{\mathcal D_0,r}/\mathcal D_\infty,\beta)\right]^{\on{vir}},
\end{align}
$\tau^\prime$ is the forgetful map
\[
\tau^\prime: \overline{M}_{0, \vec k}(\mathcal Y_{\mathcal D_0,r}/\mathcal D_\infty,\beta) \rightarrow \overline{M}_{0, m}(\mathcal D,\beta).
\]
We apply the virtual localization formula to $(\tau^\prime)_*\left[\overline{M}_{0, \vec k}(\mathcal Y_{\mathcal D_0,r}/\mathcal D_\infty,\beta)\right]^{\on{vir}}_{\mathbb C^\times}$ with respect to the $\mathbb{C}^*$-action that scales the fibers of $\mathcal Y_{\mathcal D_0,r}\to \mathcal D$. This is described explicitly in \cite{TY20c}*{Section 2.2.2} for genus $g$ invariants. We refer to \cite{TY20c}*{Section 2.2.2} for the precise localization formula and specialize it to the genus zero case. The edge contribution is trivial when $r$ is sufficiently large. Let $t$ be the equivariant parameter. Then the $t^{-1}$-coefficient of the localization formula is zero. The $t^{-1}$-coefficient the localization graphs are either one stable vertex over $0$ or one stable vertex over $\infty$. Therefore, the $t^{-1}$-coefficient is precisely:
\[
r(\tau_1)_*\left( [\overline{M}_{0,\vec k}(\sqrt[r]{L/\mathcal D},\beta)]^{\on{vir}}\right)-(\tau_2)_*\left([\overline{M}_{0,\vec k}^\sim(\mathcal D,\beta)]^{\on{vir}}\right)=0.
\]

    Furthermore, by \cite{TT}, 
\[
r(\tau_1)_*\left( [\overline{M}_{0,\vec k}(\sqrt[r]{L/\mathcal D},\beta)]^{\on{vir}}\right)=[\overline{M}_{0,m}(\mathcal D,\beta)]^{\on{vir}}
\]
\end{proof}


Now we consider the genus zero Gromov--Witten theory of the root stack $\mathcal X_{\mathcal D,r}$ with one large age (or, equivalently, genus zero relative Gromov--Witten theory of the orbifold pair $(\mathcal X,\mathcal D)$ with one negative contact). We assume that $r$ is sufficiently large and consider the degeneration formula.
\begin{proposition}\label{prop-1-neg}
Let 
\begin{align*}
   \alpha_i\in H^*_{\on{CR}}(\mathcal X_{\mathcal D,r}), \text{ and } a_i\in \mathbb Z_{\geq 0}  \text{ for } i\in \{1,\ldots, m\}. 
\end{align*}

Without loss of generality, we assume that $\alpha_1$ is the insertion that corresponds to the unique negative contact marking. Then the relative invariant with one negative contact order can be written as
\begin{align}\label{degen-formula-1-neg}
    \langle \prod_{i=2}^m \tau_{a_i}(\alpha_i)  \rangle_{\Gamma}^{(\mathcal X,\mathcal D)}=\sum \frac{\prod_{e\in E}d_e}{|\Aut(E)|} \langle  \prod_{i\in S_{1}}\tau_{a_i}(\alpha_i), \eta, \tau_{a_1}(\alpha_1)\rangle^\sim_{\Gamma^0}\langle \check{\eta}, \prod_{i\in S_{2}}\tau_{a_i}(\alpha_i) \rangle_{\Gamma^\infty}^{\bullet, (\mathcal X,\mathcal D)},
\end{align}
where $\Aut(E)$ is the permutation group of the set $\{d_1,\ldots, d_{|E|}\}$; $\langle \cdots\rangle^\sim_{\Gamma^0}$ is the rubber invariant; $\check{\eta}$ is defined by taking the Poincar\'e dual of the orbifold cohomology weights of the cohomology weighted partition $\eta$;  the second sum is over all splittings of 
\[
\{2,\ldots, m\}=S_{1}\sqcup S_{2}
\]
and all intermediate orbifold cohomology weighted partitions $\eta$.    
\end{proposition}

\begin{proof}
    We consider the degeneration to the normal cone of $\mathcal D$ in $\mathcal X$ and apply the degeneration formula for the orbifold Gromov--Witten theory of $\mathcal X_{\mathcal D,r}$ for $r$ sufficiently large. The degeneration formula expresses genus zero orbifold Gromov--Witten invariants of $\mathcal X_{\mathcal D,r}$ (equals to genus zero relative Gromov--Witten invariants of $(\mathcal X,\mathcal D)$) in terms of (possibly disconnected) invariants of $(\mathcal X,\mathcal D)$ and $(\mathcal Y_{\mathcal D_0, r}, \mathcal D_\infty)$. To obtain the formula (\ref{degen-formula-1-neg}), we just need to compute the (possibly disconnected) invariants of $(\mathcal Y_{\mathcal D_0, r}, \mathcal D_\infty)$ with orbifold ages along $\mathcal D_0$ is given by $k_i$, for $i\in S_1$ and contact orders along $\mathcal D_\infty$ is given by the weights $\eta_i$.

    Similar to the proof of Lemma \ref{lemma-gerbe-rubber}, we apply the $\mathbb C^\times$-virtual localization formula to the invariants of $(\mathcal Y_{\mathcal D_0, r}, \mathcal D_\infty)$. The difference is that we now consider invariants with one large age. In the localization graph, there is a vertex $V_0$ over $0$ that carries this unique marking with large age. We need to extract the coefficient of $t^0$, where $t$ is the equivariant parameter. All the vertex contributions have only coefficients of $t^i$, for $i<0$, except for the vertex contribution at $V_0$. For the vertex contribution at $V_0$, the highest power of $t$ is $t^0$. Therefore, there is only one stable vertex ($V_0$) in the localization graph. The $t^0$-coefficient is 
    \[
    [\overline{M}_{0,\vec l}(\sqrt[r]{L/\mathcal D},\beta)]^{\on{vir}},
    \]
    where the positive part of $\vec l$ is given by $\{k_i\}_{i\in S_1\setminus\{1\}}$ and the negative part of $\vec l$ is given by $k_1$ and $\{-\eta_i\}$. By Lemma \ref{lemma-gerbe-rubber}, this is precisely the rubber invariant in (\ref{degen-formula-1-neg}).
\end{proof}

When the root stack $\mathcal X_{\mathcal D,r}$ is a multi-root stack $X_{D,\vec r}$, we have a formula for the orbifold invariant of the snc pair $(X,D)$ with one negative contact order along certain irreducible component $D_i$ of $D$.



\section{Orbifold invariants with mid-ages}\label{sec:mid-age}

In this section, we introduce a new type of orbifold invariants which we call genus zero invariants of multi-root stacks with mid-ages, or mid-age invariants for short. Mid-age invariants were first studied in \cite{You21} when $D$ is smooth. We first recall the set-up in \cite{You21}, then we will consider more general mid-age invariants when $D=D_1+\cdots+D_n$ is snc. 

\subsection{Mid-age invariants for smooth pairs}
Let $X$ be a smooth projective variety and $D\subset X$ be a smooth divisor. Let $\vec k=(k_1,\ldots,k_m)$ be a vector of $m$ integers which satisfy 
\[
\sum_{i=1}^m k_i=\int_\beta[D].
\] 
In \cite{You21}, we considered orbifold invariants of the $r$-th root stack $\mathcal X:=X_{D,r}$ whose marking $p_i$ is mapped to a sector of the inertia stack $IX_{D,r}$ with age
\begin{align}\label{ages}
age_{p_i}(N_{\mathcal D/\mathcal X})=\left\{\begin{array}{cc}
     \frac{k_i}{r}  & \text{ if } k_i\geq 0\text{ and } i\leq m, \\
      1+\frac{k_i}{r} & \text{ if } k_i<0 \text{ and } i\leq m,\\
      k_a/r & i=m+1,\\
      k_b/r & i=m+2,
\end{array}\right.
\end{align}
where 
$k_a,k_b\in \{1,2,\ldots, r-1\}$ and 
\[
k_a+k_b=r.
\]
The ages $k_a/r$ and $k_b/r$ are called mid-ages if they satisfy the assumption of \cite{You21}*{Lemma 2.1}. This means that, for $r$ sufficiently large, we also choose $k_a$ and $k_b$ to be sufficiently large with respect to $D\cdot\beta$. For example, we may require $k_a/r$ and $k_b/r$ to be sufficiently close to $1/2$.

Following the standard approach of studying orbifold Gromov--Witten theory of root stacks in \cite{TY18}, \cite{FWY} and \cite{You21}, we can consider the degeneration to the normal cone to $D$ and apply the degeneration formula. The degeneration formula expresses the invariants of $X_{D,r}$ in terms of invariants of $(X,D)$ and $(Y_{D_0,r},D_\infty)$, where $Y:=\mathbb P(N_{D}\oplus \mathcal O)$, $D_0$ and $D_\infty$ are zero and infinity divisors of $Y$. Then we apply the virtual localization formula to study the orbifold-relative Gromov--Witten theory of $(Y_{D_0,r},D_\infty)$. A key feature of the mid-age invariants is that the pair of mid-age markings will always be in the same rubber. That is, the irreducible components containing mid-age markings are mapped into $D$. 


Let $\overline{M}_{0,\vec k,k_a,k_b,\beta}(X_{D,r})$ be the moduli space of genus $0$, degree $\beta\in H_2(X,\mathbb Q)$, orbifold stable maps to $X_{D,r}$ with ages of the markings are given by (\ref{ages}).  Let 
 \[
 m_+:=\#\{i: k_i>0\} \quad \text{and,} \quad m_-:=\#\{i: k_i<0\}.
 \]
 We consider the forgetful map
\[
\tau: \overline{M}_{0,\vec k, k_a,k_b,\beta}(X_{D,r}) \rightarrow \overline{M}_{0,m+2,\beta}(X)\times_{X^{m_++m_-+2}}D^{m_++m_-+2},
\]
which forgets the orbifold structures.  One of the main results of \cite{You21} is the following. 

 \begin{theorem}[=\cite{You21}*{Theorem 1.2}]\label{thm-mid-age}
 The genus zero cycle class
\[
r^{m_-+1}\tau_*\left(\left[\overline{M}_{0,\vec k, k_a,k_b,\beta}(X_{D,r})\right]^{\on{vir}}\right)\in A_*\left(\overline{M}_{0,m+2,\beta}(X)\times_{X^{m_++m_-+2}}D^{m_++m_-+2}\right)
\] 
of the root stack $X_{D,r}$ is constant in $r$ for sufficiently large $r$. Furthermore, there exists a positive integer $d_0$ such that this cycle class does not depend on the pair $(k_a,k_b)$ as long as $d_0<k_a\leq k_b$.
\end{theorem}

 Now we consider a slightly more general set-up. Let $\vec k=(k_1,\ldots,k_m)$ be a vector of $m$ integers which satisfy 
\[
\sum_{i=1}^m k_i=\int_\beta[D]-c
\] 
and
\[
k_a+k_b=r+c,
\]
where $c\in \mathbb Z$.
We consider the genus zero orbifold Gromov--Witten theory of $X_{D,r}$ with $(m+2)$-markings and the corresponding ages are given by (\ref{ages}). It is straightforward to see that Theorem \ref{thm-mid-age} also holds in this setting.
 \begin{theorem}\label{thm-mid-age-c}
 For $r\gg 1$, we consider a pair of positive integers $\{k_a,k_b\}$ that satisfies 
 \[
k_a+k_b=r+c,
\]
where $c\in \mathbb Z$.
The genus zero cycle class
\[
r^{m_-+1}\tau_*\left(\left[\overline{M}_{0,\vec k, k_a,k_b,\beta}(X_{D,r})\right]^{\on{vir}}\right)\in A_*\left(\overline{M}_{0,m+2,\beta}(X)\times_{X^{m_++m_-+2}}D^{m_++m_-+2}\right)
\] 
of the root stack $X_{D,r}$ is constant in $r$ for sufficiently large $r$. Furthermore, there exists a positive integer $d_0$ such that this cycle class does not depend on the pair $(k_a,k_b)$ as long as $d_0<k_a\leq k_b$.
\end{theorem}
\begin{proof}
    The proof is identical to the proof of Theorem \ref{thm-mid-age} in \cite{You21}. We briefly repeat the idea. We consider the degeneration to the normal cone of $D$ in $X$ and apply the degeneration formula. The Gromov--Witten class of $X_{D,r}$ is expressed in terms of Gromov--Witten classes of $(X,D)$ and $(Y_{D_0,r},D_\infty)$. The $r$-dependency is determined by the Gromov--Witten classes of $(Y_{D_0,r},D_\infty)$ which can be computed by $\mathbb C^\times$-virtual localization formula. As in \cite{You21}*{Lemma 2.1}, one can find a positive integer $d_0$ such that the pair of markings with mid-ages have to be distributed to the same vertex $v_0$ of the localization graph. For $m_-=0$, the fixed locus corresponding to $v_0$ is
    \[
    r[\overline{M}(\mathcal D_0)]^{\on{vir}},
    \]
    moduli space of stable maps to $\mathcal D_0$, a $\mu_r$ gerbe over $D$. By \cite{AJT15}, the pushforward of $r[\overline{M}(\mathcal D_0)]^{\on{vir}}$ is a moduli space of stable maps to $D$ which is independent of $k_a,k_b$. For $m_->0$, the idea is the same. The details are presented in \cite{You21}*{Section 2.2}. One just needs to notice that the proof of \cite{You21}*{Lemma 2.4} only requires $d_0<k_a,k_b$ and does not depend on $c=r-k_a-k_b$. This implies that the virtual cycle associated with the vertex $v_0$ is independent of $k_a$ and $k_b$. 
    \end{proof}

Note that, although the mid-age invariants in Theorem \ref{thm-mid-age} and Theorem \ref{thm-mid-age-c} are independent of sufficiently large $r$, $k_a$ and $k_b$, but they can depend on $c$. 
\begin{definition}
    The genus zero cycle class for the orbifold Gromov--Witten theory of $(X,D)$ with a pair of mid-ages is
    \[
\left[\overline{M}_{0,\vec k, \mathbf{b},-\mathbf{b}+c,\beta}(X,D)\right]^{\on{vir}}:=\lim_{r\rightarrow \infty} r^{m_-+1}\tau_*\left(\left[\overline{M}_{0,\vec k, k_a,k_b,\beta}(X_{D,r})\right]^{\on{vir}}\right),
\]
the cycle class in Theorem \ref{thm-mid-age-c} when $r$ is sufficiently large. 
\end{definition}

We will also use the notation $\langle \cdots\rangle^{(X,D)}$ or $\langle \cdots \rangle^{X_{D,\infty}}$ to denote the orbifold invariants with mid-ages. We will specify contact orders (ages), so it will be clear whether the invariants have mid-ages or not.

\begin{remark}
The symbol $\mathbf{b}$ stands for ``big" as the pair of mid-age markings may be considered as a pair of markings, one with positive contact and the other with negative contact, such that the absolute values of the contact orders of these two markings are sufficiently large and the sum of the contact orders is a constant $c$. By considering $\mathbf b$ as a sufficiently large integer, we use $\mathbf b$ and $-\mathbf b+c$ to denote the contact orders of these two mid-age markings. We can also denote the contact orders as $\mathbf b+c$ and $-\mathbf b$ because the invariant is independent of the actual value of $\mathbf b$ when it is sufficiently large.
\end{remark}

For a smooth pair $(X,D)$ there is only one broken line type. The theta functions are defined in \cite{You22} in terms of generating functions of two-point relative Gromov--Witten invariants. In this section, we will show that they can also be written in terms of generating functions of mid-age invariants.

\begin{proposition}\label{prop-mid-age-c}
  Let $(X,D)$ be a smooth log Calabi--Yau pair, the following three-point orbifold invariants with a pair of mid-ages are related to two-point relative invariants without mid-ages:
    \[
    \langle [1]_p,[1]_{\mathbf{b}},[\on{pt}]_{-\mathbf{b}+c}\rangle_{0,3,\beta}^{(X,D)}=c\langle [1]_p,[\on{pt}]_c\rangle_{0,2,\beta}^{(X,D)},
    \]
    where $p,c\in \mathbb Z_{>0}$ and $p+c=D\cdot\beta$.   
\end{proposition}
\begin{proof}
    We consider the degeneration to the normal cone of $D$ and apply the degeneration formula. Then, we apply the virtual localization formula as in \cite{You21}*{Lemma 2.2}. This is identical to the computation of the relative invariant:
\[
    \langle [1]_p,[1]_{b},[\on{pt}]_{-b+c}\rangle_{0,3,\beta}^{(X,D)},
    \]
when $b\gg D\cdot\beta$. Following \cite{You22}*{Proposition 3.3}, we have
\[
    \langle [1]_p,[1]_{\mathbf{b}},[\on{pt}]_{-\mathbf{b}+c}\rangle_{0,3,\beta}^{(X,D)}=c\langle [1]_p,[\on{pt}]_c\rangle_{0,2,\beta}^{(X,D)}.
    \]
\end{proof}

\begin{proposition}\label{prop-mid-age-sm-neg}
  Let $(X,D)$ be a smooth log Calabi--Yau pair, the following three-point orbifold invariants with mid-ages satisfy:
    \[
    \langle [1]_p,[1]_{\mathbf{b}},[\on{pt}]_{-\mathbf{b}+c}\rangle_{0,3,\beta}^{(X,D)}=\left\{\begin{array}{cc}
       0  & p\neq -c, \\
        1 & p=-c,
    \end{array}\right.
    \]
    where $p\in \mathbb Z_{>0}$, $c\in \mathbb Z_{\leq 0}$ and $p+c=D\cdot\beta$.   
\end{proposition}

\begin{proof}
    We follow the same proof as Proposition \ref{prop-mid-age-c}. The invariant is identical to the relative invariant:
\[
    \langle [1]_p,[1]_{b},[\on{pt}]_{-b+c}\rangle_{0,3,\beta}^{(X,D)},
    \]
when $b\gg D\cdot\beta$. Now we have $b-c>b,p$. Following \cite{You22}*{Proposition 3.3}, the invariant is zero unless $p=-c$. When $p=-c$, we have
\[
    \langle [1]_p,[1]_{b},[\on{pt}]_{-b-p}\rangle_{0,3,\beta}^{(X,D)}=1
\]
\end{proof}

    Therefore, \cite{You22}*{Proposition 3.3 \& Definition 3.6} can be restated in terms of mid-age invariants.

    \begin{corollary}\label{cor-theta-sm}
        The structure constants $N_{p,q,-r}^\beta$ for a smooth log Calabi--Yau pair $(X,D)$ satisfy the identity:
    \[
    N_{p,q,-r}^\beta=    \langle [1]_p,[1]_{\mathbf{b}},[\on{pt}]_{-\mathbf{b}+q-r}\rangle_{0,3,\beta}^{(X,D)}+    \langle [1]_{q},[1]_{\mathbf{b}},[\on{pt}]_{-\mathbf{b}+p-r}\rangle_{0,3,\beta}^{(X,D)}.
    \] 
For a positive integer $p$,    the theta function $\vartheta_p$ for the smooth log Calabi--Yau pair is 
    \begin{align}\label{theta-func-def-1}
\vartheta_p=x^{-p}+\sum_{n=1}^{\infty}\sum_{\beta:D\cdot\beta=n+p}n\langle [1]_p,[1]_{\mathbf{b}},[\on{pt}]_{-\mathbf{b}+n}\rangle_{0,3,\beta}^{(X,D)}t^{n+p}x^n.
\end{align}
    \end{corollary}

\subsection{Mid-age invariants for orbifold pairs}
 
We extend Theorem \ref{thm-mid-age-c} to a smooth orbifold pair $(\mathcal X,\mathcal D)$, where $\mathcal X$ is a smooth proper Deligne--Mumford stack over $\mathbb C$ with a projective coarse moduli space and $\mathcal D$ is a smooth irreducible divisor of $\mathcal X$. Let $r$ be a positive integer. We assume that $r$ is coprime with the order of any stabilizer of $\mathcal X$. Then the $r$-th root stack along $\mathcal D$,
\[
\mathcal X_{\mathcal D,r},
\]
is a smooth Deligne--Mumford stack. 

Let 
\[
\vec k=(k_1,\ldots,k_m,k_a,k_b)\in \mathbb Q^{m+2} 
\]
such that
\[
\sum_{i=1}^m k_i=\int_\beta[D]-c
\] 
and
\[
k_a+k_b=r+c,
\]
where $c\in \mathbb Q$. 
Let $\overline{M}_{0,\vec k,k_a,k_b,\beta}(\mathcal X_{\mathcal D,r})$ be the moduli space of genus $0$, degree $\beta\in H_2(\mathcal X,\mathbb Q)$, orbifold stable maps to $\mathcal X_{\mathcal D,r}$ with ages of the markings are given by
\begin{align}\label{ages-orb}
age_{p_i}(N_{\mathcal D/\mathcal X_{\mathcal D,r}})=\left\{\begin{array}{cc}
     \frac{k_i}{r}  & \text{ if } k_i\geq 0\text{ and } i\leq m; \\
      1+\frac{k_i}{r} & \text{ if } k_i<0 \text{ and } i\leq m;\\
      k_a/r & i=m+1;\\
      k_b/r & i=m+2.
\end{array}\right.
\end{align}
  Similar to the case of smooth pairs, we write
 \[
 m_+:=\#\{i: k_i>0\} \quad \text{and}\quad m_-:=\#\{i: k_i<0\}.
 \]
 We consider the forgetful map
\[
\tau: \overline{M}_{0,\vec k, k_a,k_b,\beta}(\mathcal X_{\mathcal D,r}) \rightarrow \overline{M}_{0,m+2,\beta}(\mathcal X)\times_{\mathcal X^{m_++m_-+2}}\mathcal D^{m_++m_-+2}.
\]

The following theorem is a generalization of Theorem \ref{thm-mid-age-c} to orbifold pairs. The proof is similar to the proof of Theorem \ref{thm-mid-age-c} using the orbifold degeneration formula and the virtual localization formula which are explained in details in \cite{TY20c}*{Section 2} and Theorem \ref{thm-mid-age-c}.
\begin{theorem}\label{thm-mid-age-orb}
 The genus zero cycle class
\[
r^{m_-+1}\tau_*\left(\left[\overline{M}_{0,\vec k, k_a,k_b,\beta}(\mathcal X_{\mathcal D,r})\right]^{\on{vir}}\right)\in A_*\left(\overline{M}_{0,m+2,\beta}(\mathcal X)\times_{\mathcal X^{m_++m_-+2}}\mathcal D^{m_++m_-+2}\right)
\] 
of the root stack $\mathcal X_{\mathcal D,r}$ is constant in $r$ for sufficiently large $r$. Furthermore, there exists a positive rational number $d_0$ such that this cycle class does not depend on the pair $(k_a,k_b)$ as long as $d_0<k_a\leq k_b$.
\end{theorem}


\subsection{Mid-age invariants for snc pairs}
Now we are ready to consider the Gromov--Witten theory of multi-root stacks with mid-ages. Let $X$ be a smooth projective variety over $\mathbb C$. Let
\[
D:=D_1+\cdots+D_n
\]
be an snc divisor of $X$ and $D_1,\ldots, D_n$ are smooth irreducible divisors of $X$. Let $\vec r=(r_1,\ldots, r_n)\in (\mathbb Z_{>0})^n$. We assume that $r_1,\ldots,r_n$ are pairwise coprime. The multi-root stack $X_{D,\vec r}$ is a smooth Deligne--Mumford stack. We consider the moduli space $\overline{M}_{0,\{\vec k^j\}_{j=1}^m, \vec k_a,\vec k_b,\beta}(X_{D,\vec r})$ of genus zero orbifold stable maps to $X_{D,\vec r}$ with degree $\beta\in H_2(X)$ and $(m+2)$-markings where ages of the first $m$ markings are determined by $\{\vec k^j\}_{j=1}^m$ as in Section \ref{sec:GW-multi-root}. The last two markings are mid-ages along some of $D_i$ and have large or small ages (including age $0$) along other $D_i$. By reordering $D_i$, we may assume that the last two markings are of mid-ages along $D_1,\ldots,D_{n^\prime}$, for some $n^\prime\leq n$. We also recall the following notation
\[
m_{i,-}:=\#\left\{j\in\{1,\ldots,m+2\}: k_i^j<0\right\}, \text{ for } i=1,2,\ldots, n.
\]

For each $i\in\{1,\ldots, n\}$, the root stack $X_{D,\vec r}$ can be viewed as 
\[
\left(X_{(D_1,r_1),\ldots,\widehat{(D_i,r_i)},\ldots, (D_n,r_n)}\right)_{D_i,r_i}
\]
the $r_i$-th root stack of $X_{(D_1,r_1),\ldots,\widehat{(D_i,r_i)},\ldots, (D_n,r_n)}$ along $D_i$. Therefore, applying Theorem \ref{thm-mid-age-orb} to $X_{D,\vec r}$ for each $i$, we obtain the following.
\begin{thm}\label{thm-mid-age-mul}
 The genus zero cycle class

\[
\left(\prod_{i=1}^nr_i^{m_{i,-}}\prod_{i=1}^{n^\prime}r_i\right)\tau_*\left(\left[\overline{M}_{0,\{\vec k^j\}_{j=1}^m, \vec k_a,\vec k_b,\beta}(X_{D,\vec r})\right]^{\on{vir}}\right)
\]
of the multi-root stack $X_{D,\vec r}$ is constant in $\vec r$ for sufficiently large $\vec r$. Furthermore, there exists $\vec d_0\in \mathbb Q^{n^\prime}$ such that this cycle class does not depend on the pairs $(k_{i,a}, k_{i,b})$ as long as $d_{i,0}<k_{i,a}\leq k_{i,b}$, for $i=1,\ldots, n^\prime$.
\end{thm}

\begin{definition}\label{def-mid-age}
    The genus zero cycle class for the orbifold Gromov--Witten theory of $(X,D)$ with a pair of mid-ages is
    \[
\left[\overline{M}_{0,\{\vec k^j\}_{j=1}^m, \vec{\mathbf{b}},-\vec{\mathbf{b}}+\vec c,\beta}(X,D)\right]^{\on{vir}}:=\lim_{\vec r\rightarrow \infty} \left(\prod_{i=1}^nr_i^{m_{i,-}}\prod_{i=1}^{n^\prime}r_i\right)\tau_*\left(\left[\overline{M}_{0,\{\vec k^j\}_{j=1}^m, \vec{k}_a,\vec{k}_b,\beta}(X_{D,\vec r})\right]^{\on{vir}}\right),
\]
the cycle class in Theorem \ref{thm-mid-age-mul} when $\vec r$ is sufficiently large. 
\end{definition}

Again, we will use the notation $\langle \cdots\rangle^{(X,D)}$ or $\langle \cdots \rangle^{X_{D,\infty}}$ to denote the orbifold invariants of snc pairs with mid-ages.

\begin{remark}
    The genus zero cycle class for the mid-age invariants in Definition \ref{def-mid-age} only depend on the discrete data $\{\vec k^{j}\}_{j=1}^m$, $\vec c$, $\beta$ and $k_{i,a}$, $k_{i,b}$, for $i>n^\prime$, but do not depend on $\vec r$, $ k_{i,a}$, $ k_{i,b}$, for $1\leq i \leq n^\prime$. Therefore, this cycle class is indeed intrinsic to the relative geometry of $(X,D)$. One can also define the cycle class directly using the moduli spaces of relative spaces and rubber moduli spaces without referring to orbifold Gromov--Witten theory. The advantage of considering orbifold Gromov--Witten theory is that one can apply the structural properties of orbifold Gromov--Witten theory directly. In particular, we will use the WDVV equation, the topological recursion relation and the orbifold degeneration formula.
\end{remark}

\subsection{Some properties of mid-age invariants}\label{sec:properties-mid-age}

Mid-age invariants are orbifold invariants, hence they satisfy the universal equations for orbifold Gromov--Witten invariants: string equation, divisor equation, dilaton equation, topological recursion relation (TRR), and Witten-Dijkgraaf-Verlinde-Verlinde (WDVV) equation. We state the WDVV equation and TRR which will be used in later sections. 

Let $\{\tilde{T}_{\vec s,l}\}_{\vec s, l}$ and $\{\tilde{T}_{-\vec s}^l\}_{\vec s, l}$ be dual basis of $H^*_{\on{CR}}(X_{D,\vec r})$. We use the notation $\langle \cdots \rangle^{(X,D)}$ instead of $\langle \cdots\rangle^{X_{D,\vec r}}$ to denote the invariants of $X_{D,\vec r}$ as we will always assume that $\vec r$ is sufficiently large. The WDVV equation is

\begin{align}\label{wdvv}
    &\sum \left\langle [\gamma_1]_{\vec k^1}\bar{\psi}^{a_1},[\gamma_2]_{\vec k^2}\bar{\psi}^{a_2},\prod_{j\in S_1}[\gamma_j]_{\vec k^j}\bar{\psi}^{a_j}, \tilde{T}_{\vec s,l} \right\rangle_{0,\{\vec k^j\}_{j\in S_1\cup\{1,2\}},\vec s,\beta_1}^{(X,D)}\\
\notag & \qquad \cdot\left\langle \tilde{T}_{-\vec s}^l,[\gamma_3]_{\vec k^3}\bar{\psi}^{a_3},[\gamma_4]_{\vec k^4}\bar{\psi}^{a_4},\prod_{j\in S_2}[\gamma_j]_{\vec k^j}\bar{\psi}^{a_j} \right\rangle_{0,-\vec s,\{\vec k^j\}_{j\in S_2\cup\{3,4\}},\beta_2}^{(X,D)}\\
\notag=&\sum \left\langle [\gamma_1]_{\vec k^1}\bar{\psi}^{a_1},[\gamma_3]_{\vec k^3}\bar{\psi}^{a_3}\prod_{j\in S_1}[\gamma_j]_{\vec k^j}\bar{\psi}^{a_j}, \tilde{T}_{\vec s,l} \right\rangle_{0,\{\vec k^j\}_{j\in S_1\cup\{1,3\}},\vec s,\beta_1}^{(X,D)}\\
\notag & \qquad \cdot\left\langle \tilde{T}_{-\vec s}^l,[\gamma_2]_{\vec k^2}\bar{\psi}^{a_2},[\gamma_4]_{\vec k^4}\bar{\psi}^{a_4},\prod_{j\in S_2}[\gamma_j]_{\vec k^j}\bar{\psi}^{a_j} \right\rangle_{0,-\vec s,\{\vec k^j\}_{j\in S_2\cup\{2,4\}},\beta_2}^{(X,D)},
\end{align}
where each sum is over all splittings of $\beta_1+\beta_2=\beta$, all indices $\vec s, l$ of basis, and all splittings of disjoint sets $S_1$, $S_2$ with $S_1\cup S_2=\{5,\ldots,m\}$. all of $\{\vec k^{j}\}_{j=1}^m$, as well as $\vec s$, can be of mid-ages. Note that both sides are finite sums.

We consider the special case when $m=4$, the first two markings $\{p_i\}_{i=1}^2$ are not mid-ages and the last two markings $\{p_i\}_{i=3}^4$ are mid-age markings. We require the ages of $p_3$ and $p_4$ to be mid-ages with respect to all splittings of $\beta$. The insertions for the markings are:
\[
[\gamma_1]_{\vec k^1}, [\gamma_2]_{\vec k^2}, [\gamma_3]_{\vec {\mathbf b}+\vec c}, [\gamma_4]_{-\vec {\mathbf b}}.
\]
Then the WDVV equation becomes
\begin{align*}
    &\sum \left\langle [\gamma_1]_{\vec k^1},[\gamma_2]_{\vec k^2} \tilde{T}_{\vec s,l} \right\rangle_{0,\{\vec k^1,\vec k^2,\vec s\},\beta_1}^{(X,D)} \cdot\left\langle \tilde{T}_{-\vec s}^l,[\gamma_3]_{\vec {\mathbf b}+\vec c},[\gamma_4]_{-\vec {\mathbf b}}\right\rangle_{0,\{-\vec s,\vec {\mathbf b}+\vec c,\vec{\mathbf b}\},\beta_2}^{(X,D)}\\
\notag=&\sum \left\langle [\gamma_1]_{\vec k^1},[\gamma_3]_{\vec {\mathbf b}+\vec c} \tilde{T}_{\vec s,l} \right\rangle_{0,\{\vec k^1,\vec{\mathbf b}+\vec c,\vec s\},\beta_1}^{(X,D)}\cdot\left\langle \tilde{T}_{-\vec s}^l,[\gamma_2]_{\vec k^2},[\gamma_4]_{-\vec{\mathbf b}} \right\rangle_{0,\{-\vec s,\vec k^2, -\vec {\mathbf b}\},\beta_2}^{(X,D)}.
\end{align*}
Note that $\vec s$ on the first line is not mid-age and the $\vec s$ on the second line is mid-age. Therefore, we can rewrite the WDVV equation as
\begin{align}\label{wdvv-4}
    &\sum \left\langle [\gamma_1]_{\vec k^1},[\gamma_2]_{\vec k^2} \tilde{T}_{\vec s,l} \right\rangle_{0,\{\vec k^1,\vec k^2,\vec s\},\beta_1}^{(X,D)} \cdot\left\langle \tilde{T}_{-\vec s}^l,[\gamma_3]_{\vec {\mathbf b}+\vec c},[\gamma_4]_{-\vec {\mathbf b}}\right\rangle_{0,\{-\vec s,\vec {\mathbf b}+\vec c,\vec{\mathbf b}\},\beta_2}^{(X,D)}\\
\notag=&\sum \left\langle [\gamma_1]_{\vec k^1},[\gamma_3]_{\vec {\mathbf b}+\vec c} \tilde{T}_{\vec {\mathbf b}+\vec s,l} \right\rangle_{0,\{\vec k^1,\vec{\mathbf b}+\vec c,\vec{\mathbf b}+\vec s\},\beta_1}^{(X,D)}\cdot\left\langle \tilde{T}_{-\vec{\mathbf b}-\vec s}^l,[\gamma_2]_{\vec k^2},[\gamma_4]_{-\vec{\mathbf b}} \right\rangle_{0,\{-\vec{\mathbf b}-\vec s,\vec k^2, -\vec {\mathbf b}\},\beta_2}^{(X,D)},
\end{align}
where $\vec s$ on the first and the second line are not mid-ages. Each sum is over all splittings of $\beta_1+\beta_2=\beta$ and all indices $\vec s, l$ of basis.

Similarly, TRR for invariants with a pair of mid-ages can be written as follows:
\begin{equation}\label{trr}
    \begin{split}
        &\left\langle [\gamma_1]_{\vec s^1}\bar{\psi}^{a_1+1},[\gamma_2]_{\vec{\mathbf b}+\vec c}\bar{\psi}^{a_2},[\gamma_3]_{-\vec{\mathbf b}}\bar{\psi}^{a_3},\ldots, [\gamma_m]_{\vec s^m}\bar{\psi}^{a_m} \right\rangle_{0,\{\vec s^j\}_{j=1}^m,\beta}^{(X,D)}\\
=&\sum \left\langle [\gamma_1]_{\vec s^1}\bar{\psi}^{a_1},\prod_{j\in S_1}[\gamma_j]_{\vec s^j}\bar{\psi}^{a_j}, \tilde{T}_{\vec s,k} \right\rangle_{0,\{\vec s^j\}_{j\in S_1\cup\{1\}},\vec s,\beta_1}^{(X,D)}\\
& \qquad \cdot\left\langle \tilde{T}_{-\vec s}^k,[\gamma_2]_{\vec {\mathbf b}+\vec c}\bar{\psi}^{a_2},[\gamma_3]_{\vec {\mathbf b}}\bar{\psi}^{a_3},\prod_{j\in S_2}[\gamma_j]_{\vec s^j}\bar{\psi}^{a_j} \right\rangle_{0,-\vec s,\{\vec s^j\}_{j\in S_2\cup\{2,3\}},\beta_2}^{(X,D)},
    \end{split}
\end{equation}
where the sum is over all splittings of $\beta_1+\beta_2=\beta$, all indices $\vec s, k$ of basis, and all splittings of disjoint sets $S_1$, $S_2$ with $S_1\cup S_2=\{4,\ldots,m\}$. Note that $\vec s$ are not mid-ages.

We also would like to point out that the degeneration formula for multi-root stacks mentioned in Section \ref{sec:deg} also holds for mid-age invariants. We need to pay attention to the contact orders and ages of the orbifold cohomology weighted partition $\eta$. Recall that if we consider the degeneration to the normal cone of $D_1$ and the $i$-th marking of $\eta$ has contact orders and ages
\[
(\eta_{1}^i,\eta_{2}^i,\ldots, \eta_{n}^i),
\]
then the $i$-th marking of $\eta^\vee$ has contact orders and ages 
\[
(\eta_{1}^i,-\eta_{2}^i,\ldots, -\eta_{n}^i).
\]
If the original invariant has mid-ages, then, for $j\geq 2$, it is possible that $\eta_{j}^i$ is mid-age. Note that $D_1$ is the relative divisor, there are no mid-ages with respect to $D_1$ in $\eta$. 


For the rest of this section, we prove some identities of mid-age invariants which are essential for the theta functions to be well-defined. For this purpose, we assume that $(X,D)$ is an snc log Calabi--Yau pair. However, lots of the arguments below work without the log Calabi--Yau assumption.

Let $\vec c=(c_1,\ldots, c_n)\in \mathbb Z^n$ with $D_{\vec c}\neq \emptyset$. One can consider orbifold invariants with a pair of mid-ages along $D_{i}$, where $c_i\neq 0$, such that the mid-ages $\vec k_a$ and $\vec k_b$ satisfies
\[
k_{a_i}+k_{b_i}=r_i+c_i.
\]
We denote this pair of mid-ages by $\vec {\mathbf b}$ and $-\vec{\mathbf b}-\vec c$. 

We would like to encode the information that there are irreducible components of the curves mapping into a lowest dimensional stratum of $D$ that are included in $D_{\vec c}$. Therefore, we would like to have mid-ages along this lowest dimensional stratum of $D$. If $D_{\vec c}$ is already a lowest dimensional stratum of $D$, then the choice of the mid-age $\vec{\mathbf b}$ is unique.

 If $D_{\vec c}$ is not a lowest dimensional stratum of $D$, we can also consider orbifold invariants with mid-ages along a set of $\{D_i\}_{i\in S_{\vec c}}$ such that $\cap_{i\in S_{\vec c}} D_{i}$ is a lowest dimensional (nonempty) stratum inside $D_{\vec c}$. We denote this pair of mid-ages by $\vec {\mathbf b}$ and $-\vec{\mathbf b}-\vec c$. 
 Suppose there is another set of $\{D_i\}_{i\in S^\prime_{\vec c}}$ such that $\cap_{i\in S^\prime_{\vec c}} D_{i}$ is another lowest dimensional (nonempty) stratum inside $D_{\vec c}$, we denote this pair of mid-ages by $\vec {\mathbf b}^\prime$ and $-\vec{\mathbf b}^\prime-\vec c$. We claim the mid-age invariant does not depend on the choice of the lowest dimensional strata.

  Without loss of generality, we assume that 
    \[
    \vec c=(c_1,\ldots, c_l,0,\ldots,0),
    \]
    where $c_i\neq 0$ for $i\in \{1,\ldots, l\}$, and
    \[
    \vec{\mathbf b}=(\mathbf b_1,\ldots, \mathbf b_s,0,\ldots, 0),
    \]
    where $\mathbf b_i$, for $i\in\{1,\ldots,s\}$, stands for mid-age along $D_1,\ldots, D_s$. Note that $s > l$ and $\cap_{i=1}^sD_s$ is a lowest dimensional stratum. We can assume that there is another lowest dimensional stratum  corresponding to a vector of mid-ages of the form:
    \[
    \vec{\mathbf b}^\prime=(\mathbf b_1,\ldots, \mathbf b_{s-1},0,\mathbf b_{s+1},0,\ldots, 0).
    \]
    
\begin{thm}\label{thm-indep-mid}
    The following identity holds for mid-age invariants
    \begin{align}\label{iden-indep-mid}
    \langle [1]_{\vec p},[1]_{\vec{\mathbf{b}}},[\on{pt}]_{-\vec{\mathbf{b}}+\vec c}\rangle_{0,3,\beta}^{(X,D)}= \langle [1]_{\vec p},[1]_{\vec{\mathbf{b}}^\prime},[\on{pt}]_{-\vec{\mathbf{b}}^\prime+\vec c}\rangle_{0,3,\beta}^{(X,D)}.
    \end{align}
\end{thm}

\begin{proof}

    We would like to apply the degeneration to the normal cone of $D_i$ in $X$ and compare the degeneration formula for LHS and RHS of (\ref{iden-indep-mid}).
\if{
    First, we consider the degeneration to the normal cone of $D_1$ and apply the degeneration formula. For the invariant in (\ref{iden-indep-mid}), the degeneration formula for mid-age invariants  requires that the second marking $p_2$ and the third marking $p_3$ are in the same rubber component over $D_1$. The possible difference between the LHS and the RHS of (\ref{iden-indep-mid}) is the invariant associated with this rubber component. 
    The rubber invariant containing these two markings that appears in the degeneration formula is
    \[
    \langle \eta,[\on{pt}]_{(\mathbf{b}_1-c_1,-\widehat{\vec{\mathbf{b}}}_1+\widehat{\vec c}_1 )} | \, | [1]_{\vec{\mathbf{b}}}\rangle^{\sim, (\mathbb P_{\mathcal D_1},\mathcal D_{1,0}+\mathcal D_{1,\infty})}
    \]
    or 
     \[
    \langle \eta, [\on{pt}]_{(\mathbf{b}_1-c_1,-\widehat{\vec{\mathbf{b}}}_1+\widehat{\vec c}_1 )} | \, | [1]_{\vec p},  [1]_{\vec{\mathbf{b}}}\rangle^{\sim, (\mathbb P_{\mathcal D_1},\mathcal D_{1,0}+\mathcal D_{1,\infty})},
    \]
    where......

These contributions are precisely the difference in the degeneration formulas for the invariants in (\ref{iden-indep-mid}). 
}\fi

Recall that $n$ is the number of irreducible components of the snc divisor $D=D_1+\ldots+D_n$. It is trivial when $n=2$. We start with the case when $n=3$ such that $D_{12}:=D_1\cap D_2\neq \emptyset$, $D_{13}:=D_1\cap D_3 \neq \emptyset$ and $D_{123}:=D_1\cap D_2\cap D_3=\emptyset$. Then we need to prove
    \begin{align}\label{iden-mid-age-3}
        \langle [1]_{(p_1,p_2,p_3)},[1]_{(\mathbf b_1+c,\mathbf b_2,0)},[\on{pt}]_{(-\mathbf b_1,-\mathbf b_2,0)}\rangle= \langle [1]_{(p_1,p_2,p_3)},[1]_{(\mathbf b_1+c,0,\mathbf b_3)},[\on{pt}]_{(-\mathbf b_1,0,-\mathbf b_3)}\rangle,
    \end{align}
    where we ignore the superscript and subscript of $\langle \cdots \rangle_{0,3,\beta}^{(X,D)}$ when there is no ambiguity. 
    
We also introduce an invariant that has mid-age only along $D_1$:
\begin{align}\label{inv-mid-age-D-1}
\langle [1]_{(p_1,p_2,p_3)},[D_2-D_3]_{(\mathbf b_1+c,0,0)},[D_1^\vee]_{(-\mathbf b_1,0,0)}\rangle,
\end{align}
where $[D_2-D_3],[D_1]\in H^2(D_1)$ are the restrictions of $[D_2-D_3], [D_1]\in H^2(X)$ under the inclusion map and $[D_1^\vee]$ is the Poincar\'e dual of $[D_1]$ in $H^*(D_1)$.


We can consider the degeneration to the normal cone of $D_1$ and apply degeneration formulas to the above three invariants in (\ref{iden-mid-age-3}) and (\ref{inv-mid-age-D-1}). The degeneration formulas for these three invariants are the same except for the rubber invariants that contain the second and the third markings. The rubber moduli can be push forwarded to the moduli of stable maps to $(D_1,D_{12}+D_{13})$ and the rubber invariants become (orbifold) invariants of $(D_1,D_{12}+D_{13})$:
\begin{align}\label{inv-mid-age-deg-1}
    \langle \eta, [1]_{(\mathbf b_2,0)},[\on{pt}]_{(-\mathbf b_2,0)}\rangle^{(D_1,D_{12}+D_{13})},
\end{align}
\begin{align}\label{inv-mid-age-deg-2}
\langle \eta, [1]_{(0,\mathbf b_3)},[\on{pt}]_{(0,-\mathbf b_3)}\rangle^{(D_1,D_{12}+D_{13})}, 
\end{align}
and
\begin{align}\label{inv-mid-age-deg-12}
\langle \eta, [D_2-D_3]_{(0,0)},[D_1^\vee]_{(0,0)}\rangle^{(D_1,D_{12}+D_{13})},    
\end{align}
respectively, where $\eta$ may have several markings with possibly nonzero contact orders to $D_{12}$ or $D_{13}$. The insertions $\eta$ also include the first markings of the original invariants if they are distributed to the rubber. Note that $D_{12}$ and $D_{13}$ are disjoint divisors. I will show in Lemma \ref{lemma-mid-age-0} that the invariants (\ref{inv-mid-age-deg-1}) (resp. (\ref{inv-mid-age-deg-2})) are zero unless they are of degree zero and the entire curve is mapped to $D_{12}$ (resp. $D_{13}$).

Now we are back to the invariant (\ref{inv-mid-age-D-1}), which can be written as
\[
\sum \on{Cont}(\cdots)\langle \eta, [D_2-D_3]_{(0,0)},[D_1^\vee]_{(0,0)}\rangle^{(D_1,D_{12}+D_{13})},
\]
where the sum is over all possible bipartite graphs in the degeneration formula and $\on{Cont}(\cdots)$ is the contribution other than the invariant of $(D_1,D_{12}+D_{13})$. Note that $\on{Cont}(\cdots)$ are the same for all invariants in (\ref{inv-mid-age-deg-1}) and (\ref{inv-mid-age-deg-2}).

We restrict to the degeneration graphs of the invariant (\ref{inv-mid-age-D-1}) when the invariants of $(D_1,D_{12}+D_{13})$ are of degree zero. The invariants are zero when the degrees are zero. On the other hand, there are nonzero contributions when the curves are entirely in $D_{12}$ or $D_{13}$. The contributions when the curves are entirely in $D_{12}$ cancel with the contributions when the curve are entirely in $D_{13}$. They, respectively, precisely correspond to the degeneration graph of the invariants on the LHS and the RHS of (\ref{iden-mid-age-3}).


    For the general case when $n>3$, recall that we assume 
\[
\vec c=(c_1,\ldots, c_{l},0\ldots,0),
\]
and    
\[
\vec{\mathbf b}=(\mathbf b_1,\ldots, \mathbf b_s,0,\ldots,0),
\]
where $l<s$, $\mathbf b_i$ stands for mid-age along $D_i$ and $\cap_{i=1}^s D_{i}$ is a lowest dimensional stratum. Without loss of generality, we assume that
\[
\vec{\mathbf b}^\prime=(\mathbf b_1,\ldots, \mathbf b_{s-1},0,\mathbf b_{s+1},0,\ldots,0),
\]
where $\cap_{i=1}^{s-1} D_{i} \cap D_{s+1}$ is another lowest dimensional stratum. Then we need to prove the following identity:
\begin{equation}\label{iden-indep-mid-s}
\begin{split}
    &\langle [1]_{\vec p},[1]_{(\mathbf b_1,\ldots, \mathbf b_s,0,\ldots,0)},[\on{pt}]_{(-\mathbf b_1+c_1,\ldots,-\mathbf b_{l}+c_{l},-\mathbf b_{l+1},\ldots, -\mathbf b_{s-1}, -\mathbf b_s,0,\ldots,0)}\rangle\\
=   &\langle [1]_{\vec p},[1]_{(\mathbf b_1,\ldots, \mathbf b_{s-1},0,\mathbf b_{s+1},0,\ldots,0)},[\on{pt}]_{(-\mathbf b_1+c_1,\ldots, -\mathbf b_{l}+c_{l},-\mathbf b_{l+1},\ldots, -\mathbf b_{s-1},0,-\mathbf b_{s+1},0,\ldots,0)}\rangle.    
\end{split}
    \end{equation}
We also introduce the invariant
\begin{align}\label{inv-mid-age-s-1}
    \langle [1]_{\vec p},[D_{s}-D_{s+1}]_{(\mathbf b_1,\ldots, \mathbf b_{s-1},0,\ldots,0)},[D_{s-1}^\vee]_{(-\mathbf b_1+c_1,\ldots,-\mathbf b_{l}+c_{l},-\mathbf b_{l+1},\ldots, -\mathbf b_{s-1},0,\ldots,0)}\rangle.
\end{align}

Similar to the case when $n=3$, we consider the degeneration to the normal cone of $D_1$ and apply the degeneration formula to invariants in (\ref{iden-indep-mid-s}) and (\ref{inv-mid-age-s-1}). The second and the third markings are distributed to the rubber. Then we pushforward the rubber moduli to the moduli of stable maps to (orbifold) $D_1$. We repeat this process by considering degeneration to the normal cone of $D_1\cap D_2$ in $D_1$ and applying the degeneration formula to the orbifold invariant of $D_1$. The pair of mid-age markings are distributed to the same rubber. After pushing forward the rubber moduli, on the LHS of (\ref{iden-indep-mid-s}), we have the invariant
\[
\langle \eta, [1]_{(\mathbf b_3+c_3,\ldots,\mathbf b_{s},0,\ldots,0)},[\on{pt}]_{(-\mathbf b_3,\ldots,-\mathbf b_{s},0,\ldots,0)}\rangle^{(D_1\cap D_2,D)},
\]
where we simply denote $D|_{D_1\cap D_2}$ by $D$.

We repeat this process to $D_{3}, \ldots, D_{s-1}$. Then the pair of mid-age markings are mapped into the rubber moduli over $\cap_{i=1}^{s-1} D_i$. After pushing forward the rubber moduli, we have the invariant
\begin{align}\label{inv-mid-age-deg-s}
\langle \eta, [1]_{(\mathbf b_{s},0,\ldots,0)},[\on{pt}]_{(-\mathbf b_{s},0,\ldots,0)}\rangle^{(\cap_{i=1}^{s-1} D_i,D)},
\end{align}  
where we denote $D|_{\cap_{i=1}^{s-1} D_i}$ by $D$ for simplicity. Similarly, we have the invariant
\begin{align}\label{inv-mid-age-deg-s1}
\langle \eta, [1]_{(0,\mathbf b_{s+1},0,\ldots,0)},[\on{pt}]_{(0,-\mathbf b_{s+1},0,\ldots,0)}\rangle^{(\cap_{i=1}^{s-1} D_i,D)},
\end{align} 
and
\begin{align}\label{inv-mid-age-deg-s-1}
\langle \eta, [D_{s}-D_{s+1}]_{(0,0,\ldots,0)},[D_{s-1}^\vee]_{(0,0,\ldots,0)}\rangle^{(\cap_{i=1}^{s-1} D_i,D)}.
\end{align} 

As the above invariants (\ref{inv-mid-age-deg-s}), (\ref{inv-mid-age-deg-s1}) and (\ref{inv-mid-age-deg-s-1}) for $(\cap_{i=1}^{s-1} D_i,D)$ may have several negative contact orders. After further degenerations and pushforwards, we get descendant invariants in general. In the language of orbifold Gromov--Witten theory, we have Hurwitz--Hodge integrals, after pushing forward, we have descendant invariants in general. Therefore, there can be descendant classes in $\eta$ for the invariant $(\cap_{i=1}^{s-1} D_i,D)$. 
Descendant classes can be removed by applying the topological recursion relation (TRR). As in (\ref{trr}), we consider the pair of mid-age markings as the second and the third markings of TRR, so these two mid-age markings will always stay together. After applying TRR to remove possible descendant classes, the rest of the argument is the same as the case when $n=3$.
    
\end{proof}


It remains to prove the following lemma.

\begin{lemma}\label{lemma-mid-age-0}
    For a log Calabi--Yau pair $(Y,E_1+E_2)$, where $E_1$ and $E_2$ are disjoint, the following genus zero mid-age invariants 
    \begin{align}\label{inv-eta-mid}
    \langle \eta, [1]_{(\mathbf b,0)},[\on{pt}]_{(-\mathbf b,0)}\rangle_{0,|\eta|+2,\beta}^{(Y,E_1+E_2)},
    \end{align}
  where $\eta$ may have several markings with possibly nonzero contact orders to $E_1$ or $E_2$, is zero unless the degree is zero. In this case, the entire curve maps into the divisor $E_1$.
\end{lemma}
\begin{proof}
    As the divisors are disjoint, we can apply the graph sum definition of relative Gromov--Witten invariants with negative contact orders in \cite{FWY}. Recall that the definition in \cite{FWY} comes from the degeneration to the normal cone to $E_1$ (or $E_2$) and the virtual localization computation. The invariant (\ref{inv-eta-mid}) equal to
    \[
    \sum_{\mathfrak G\in \mathcal B_\Gamma}\frac{\prod_{e\in E}d_e}{|\Aut(E)|} \sum_\xi \langle \eta_1,\xi\rangle_{\Gamma^\infty}^{\bullet, \mathfrak c_\Gamma, (Y,E_1+E_2)}\langle \check{\xi}, [\on{pt}]_{(\mathbf b,0)}|\eta_2,  |[1]_{(\mathbf b,0)}\rangle_{\Gamma^0}^{\sim,\mathfrak c_\Gamma},
    \]
    where $\xi$ is a cohomology weighted partition and $\langle \rangle^\sim$ is a rubber invariant over $E_1$.

    We consider the pushforward of the rubber moduli to the moduli of stable maps to $E_1$, the rubber invariant becomes
    \[
    \langle \check{\xi}, [\on{pt}]_{0},\eta_2,  [1]_{0}\rangle^{E_1},
    \]
    where there are some descendant classes in the insertion if $\eta_2$ carried negative contact orders with respect to $E_1$ from the original invariant. Let $m_-$ be the number of markings with negative contact orders (not including the mid-age markings) that are distributed to the rubber. Then the number of descendant classes is no more than $m_-$, by the graph sum definition of the relative Gromov--Witten invariants with negative contact orders in \cite{FWY}. The descendant classes can be removed by the dilaton equation or the string equation. Note that there is already a point constraint and the virtual dimension of the moduli space is
    \[
    \dim_{\mathbb C} E_1+m-3,
    \]
    where $m$ is the number of markings, therefore we can not have classes like $\psi^k$, for $k>1$ or $[\gamma]\psi$, for $\gamma\in H^*(E_1)\setminus H^0(E_1)$ because it will requires more insertions with identity insertions, then one can apply the string equations again. In the end, the invariant can be reduced to the three-point invariant
    \[
    \langle [\on{pt}],[1],[1]\rangle_{0,3,\beta}^{E_1},
    \]
    which is zero unless the degree is zero. Since the degree is zero, if $\eta_2$ includes negative contact orders, $\eta_2$ must also include at least one insertion with positive contact order. Recall that $m_-$ is the number of negative contact orders in $\eta_2$. It is also the maximal number of descendant classes that can appear. It means that we can only apply the string equation $m_-$-times, which removes $m_-$-markings. Since $\eta_2$ must have at least $m_-+1$ markings (because at least one of them is of positive contact), to end up with a three-point invariant, we need to require $\xi$ to be empty. Therefore, the entire curve has to be in $E_1$ and $\eta_2=\eta$. Note that it also means that indeed $\eta$ can not include markings with nonzero contact orders to $E_2$. 
\end{proof}


Theorem \ref{thm-indep-mid} means that for any $\vec c\in \mathbb Z^n$, there is a well-defined mid-age invariant
\[
\langle [1]_{\vec p},[1]_{\vec{\mathbf{b}}},[\on{pt}]_{-\vec{\mathbf{b}}+\vec c}\rangle_{0,3,\beta}^{(X,D)}
\]
defined by having mid-age along a set of $\{D_i\}_{i\in S_{\vec c}}$ such that $\cap_{i\in S_{\vec c}} D_{i}$ is a lowest dimensional (nonempty) stratum that contained in $D_{\vec c}$. The definition does not depend on the choice of the lowest dimensional strata. 

\begin{remark}
We can explain Theorem \ref{thm-indep-mid} in a different language. Let $\vec {|c|}=(|c_1|,\ldots, |c_n|)\in B(\mathbb Z)$, If $\vec{|c|}$ is in the interior of a maximal dimensional cone, then there is a unique way to define the mid-age invariants. If $\vec{|c|}$ is on a wall, then one can still define the mid-age invariants by choosing one of the maximal dimensional cones that are separated by the wall. The definition is independent of the choice of the maximal dimensional cone. Comparing with \cite{GS21}, this is similar to the consistency of the canonical wall structures that guarantees that the local description $\vartheta_{\vec p}(\mathrm x)$ for theta functions patch together to give a global function.   
\end{remark}


\if{
\begin{proposition}
  Let $(X,D)$ be an snc Calabi--Yau pair, the following three-point orbifold invariants with mid-ages satisfy:
    \[
    \langle [1]_{\vec p},[1]_{\vec{\mathbf{b}}},[\on{pt}]_{-\vec{\mathbf{b}}-\vec p}\rangle_{0,3,\beta}^{(X,D)}=1,
    \]
    where $\vec p\in B(\mathbb Z)$ and $D\cdot\beta=\vec 0$.   
\end{proposition}
\begin{proof}
    Without loss of generality, we assume the following 
    \begin{itemize}
        \item $\vec p=(p_1,\ldots,p_k,0,\ldots, 0)$ where $p_i>0$ for $i\in \{1,\ldots, k\}$;
        \item $\vec{\mathbf{b}}=(\mathbf{b}_1,\ldots, \mathbf{b}_l,0,\ldots, 0)$, where $l\geq k$;
        \item $D_{\vec{\mathbf{b}}}=\cap_{i=1}^lD_{i}$ is a lowest dimensional stratum of $D$.
    \end{itemize}  
    We apply the degeneration to the normal cones to $D_i$ repeatedly for $i\in\{1,\ldots, l\}$. We first apply the degeneration to the normal cones to $D_1$, the degeneration formula gives
    \begin{align*}
    \langle [1]_{\vec p},[1]_{\vec{\mathbf{b}}},[\on{pt}]_{-\vec{\mathbf{b}}-\vec p}\rangle_{0,3,\beta}^{(X,D)}
    =&\langle \eta\rangle_{0,|\eta|,\beta_1}^{(X,D)}\langle \eta^\vee, [\on{pt}]_{(\vec{\mathbf{b}}+p_1,-\widehat{\vec{\mathbf{b}}}-\widehat{\vec p})}, [1]_{\vec p},[1]_{\vec{\mathbf{b}}},\rangle_{0,3+|\eta|,\beta_2}^{(\mathbb P_{\mathcal D_1},\mathcal D_{1,0}+\mathcal D_{1,\infty})}\\
    &+\langle [1]_{\vec p},\eta\rangle_{0,1+|\eta|,\beta_1}^{(X,D)}\langle \eta^\vee, [\on{pt}]_{(\vec{\mathbf{b}}+p_1,-\widehat{\vec{\mathbf{b}}}-\widehat{\vec p})}, [1]_{\vec{\mathbf{b}}},\rangle_{0,2+|\eta|,\beta_2}^{(\mathbb P_{\mathcal D_1},\mathcal D_{1,0}+\mathcal D_{1,\infty})},
    \end{align*}
    where......
\end{proof}
}\fi

\section{Orbifold theta functions for snc pairs}

In this section, we define theta functions for snc pairs using the orbifold invariants with mid-ages introduced in Section \ref{sec:mid-age}. We call them orbifold theta functions to emphasize that they may be different from the theta functions defined by punctured logarithmic invariants in \cite{GS21}. Furthermore, we do not assume that $D$ has a zero dimensional stratum.

Recall that for a cone $\sigma$ in $B(\mathbb Z)$, we write $D_{\sigma}$ for the stratum of $D$ corresponding to $\sigma$.
\begin{definition}\label{def-theta-mid-age}
For an snc log Calabi--Yau pair $(X,D)$,    the theta functions are defined as follows: Fix $\vec p \in B(\mathbb Z)\setminus \{0\}$. Let $\sigma_{\on{max}}\in \Sigma(X)$ be a maximal cone of $\Sigma(X)$ such that $\mathrm x\in \sigma_{\on{max}}$, then
    \begin{align}\label{iden-def-theta-mid-age}
   \vartheta_{\vec p}(\mathrm x)&:=
   \sum_{\vec k\in {\mathbb Z}^n} \sum_{\beta: D_i\cdot \beta=k_i+p_i}N_{\vec p, \vec{\mathbf b}+\vec k, -\vec{\mathbf{b}}}^\beta t^{\beta} x^{\vec k},
    \end{align}
where $N_{\vec p, \vec{\mathbf b}+\vec k, -\vec{\mathbf{b}}}^\beta$ are orbifold invariants with mid-ages along the divisors $D_{\mathbf b_i}$ such that $D_{\vec {\mathbf b}}=D_{\sigma_{\on{max}}}$ is a lowest dimensional stratum in $D_{\vec k}$; $x^{\vec k}:=x_1^{k_1}\cdots x_n^{k_n}$. If $D_{\vec k}\not\supseteq D_{\sigma_{\on{max}}}$, then $N_{\vec p, \vec{\mathbf b}+\vec k, -\vec{\mathbf{b}}}^\beta=0$.
\end{definition}


Note that by Theorem \ref{thm-indep-mid}, the invariants in Definition \ref{def-theta-mid-age} are well-defined. That is, this is well-defined when $\mathrm x$ is on a wall.



\begin{remark}
    The mid-age invariant $N_{\vec p, \vec{\mathbf b}+\vec k, -\vec{\mathbf{b}}}^{\beta}$ encodes the contact order $\vec k$ and requires the second and the third markings are in the irreducible components mapping into the lowest dimensional stratum $D_{\vec {\mathbf b}}$. This may be considered as a refinement of orbifold invariants that correspond to the broken line type in \cite{GS21}.
\end{remark}

We will explain how the theta functions in Definition \ref{def-theta-mid-age} satisfies the multiplication rule:
\begin{align}\label{multi-rule-m}
\vartheta_{\vec p}\star \vartheta_{\vec q}=\sum_{\vec r\in B(\mathbb Z)} \sum_{\beta}N_{\vec p, \vec q, -\vec r}^{\beta}t^\beta\vartheta_{\vec r}.
\end{align}
By the definition of the theta functions, the LHS of (\ref{multi-rule-m}) is
\begin{align*}
    \vartheta_{\vec p}\star \vartheta_{\vec q}
    =&\left(\sum_{\vec k \in \mathbb Z^n}\sum_{\beta_1} N_{\vec p, \vec{\mathbf b}+\vec k, -\vec{\mathbf{b}}}^{\beta_1} t^{\beta_1} x^{\vec k}\right)\left(\sum_{\vec l \in \mathbb Z^n}\sum_{\beta_2} N_{\vec q, \vec{\mathbf b}+\vec l, -\vec{\mathbf{b}}}^{\beta_2}t^{\beta_2} x^{\vec l}\right)\\
=& \sum_{\vec k,\vec l\in \mathbb Z^n}\sum_{\beta_1,\beta_2}N_{\vec p, \vec{\mathbf b}+\vec k, -\vec{\mathbf{b}}}^{\beta_1}N_{\vec q, \vec{\mathbf b}+\vec l, -\vec{\mathbf{b}}}^{\beta_2}t^{\beta_1+\beta_2}x^{\vec k+\vec l}.
\end{align*}

On the other hand, the RHS of (\ref{multi-rule-m}) can be written as
\begin{align*}
\sum_{\vec r\in B(\mathbb Z)} \sum_{\beta_1}N_{\vec p, \vec q, -\vec r}^{\beta_1}t^{\beta_1}\vartheta_{\vec r}
=\sum_{\vec r\in B(\mathbb Z),\vec s\in \mathbb Z^n}\sum_{\beta_1,\beta_2} N_{\vec p, \vec q, -\vec r}^{\beta_1}N_{\vec r,\vec{\mathbf b}+\vec s,-\vec{\mathbf b}}^{\beta_2}t^{\beta_1+\beta_2} x^{\vec s}.
\end{align*}

Therefore, to equate the coefficients of $ t^{\beta}x^{\vec s}$, where $D\cdot\beta=\vec s+\vec p+\vec q$, we need to prove the following.
\begin{thm}\label{thm-wdvv}
For $D\cdot\beta=\vec s+\vec p+\vec q$ and $\vec s\in \mathbb Z^n$, we have    
\begin{align}
    \sum_{\vec k,\vec l\in \mathbb Z^n:\vec k+\vec l=\vec s}\sum_{\beta_1+\beta_2=\beta}N_{\vec p, \vec{\mathbf b}+\vec k, -\vec{\mathbf{b}}}^{\beta_1}N_{\vec q, \vec{\mathbf b}+\vec l, -\vec{\mathbf{b}}}^{\beta_2}=\sum_{\vec r\in B(\mathbb Z),\vec s\in \mathbb Z^n}\sum_{\beta_1+\beta_2=\beta} N_{\vec p, \vec q, -\vec r}^{\beta_1}N_{\vec r,\vec{\mathbf b}+\vec s,-\vec{\mathbf b}}^{\beta_2}.
\end{align}

\end{thm}
\begin{proof}
We apply the WDVV equation (\ref{wdvv-4}) for orbifold invariants with mid-ages. We consider the first four markings with insertions
\[
[1]_{\vec p}, \quad [1]_{\vec q}, \quad [1]_{\vec {\mathbf b}+\vec s}, \quad [\on{pt}]_{-\vec {\mathbf b}}.
\]    
The WDVV equation gives the following identity:
\begin{align}\label{iden-wdvv-mid-age}
    \sum \langle [1]_{\vec p},[1]_{\vec {\mathbf b}+\vec s},[\gamma]_{-\vec r}\rangle\langle [\gamma^\vee]_{\vec r}, [1]_{\vec q}, [\on{pt}]_{-\vec {\mathbf b}}\rangle
    = \sum \langle [1]_{\vec p}, [1]_{\vec q}, [\gamma]_{-\vec r}\rangle \langle [\gamma^\vee]_{\vec r}, [1]_{\vec {\mathbf b}+\vec s}, [\on{pt}]_{-\vec {\mathbf b}} \rangle,
\end{align}
where the sum is over all splittings of $\beta_1+\beta_2=\beta$, all $[\gamma]_{\vec r}$ of the basis. As pointed out in Section \ref{sec:properties-mid-age}, $\vec r$'s on the LHS are mid-ages and $\vec r$'s  on the RHS are not mid-ages.

We first consider the RHS of (\ref{iden-wdvv-mid-age}). We examine the invariant 
\[
\langle [1]_{\vec p}, [1]_{\vec q}, [\gamma]_{-\vec r}\rangle.
\]
By the virtual dimension constraint (\ref{vir-dim}), we have
\begin{align}\label{vir-dim-3}
\dim_{\mathbb C}X-3+3=0+0+d/2+\#\{i:r_i<0\},
\end{align}
where $d$ is the real degree of the cohomology class $\gamma\in \mathfrak H_{-\vec r}:=H^*(D_{-\vec r})$. Therefore, 
\[
d/2\ \leq \dim_{\mathbb C}X-\#\{i:r_i\neq 0\}.
\]
Then the RHS of (\ref{vir-dim-3}) is
\[
d/2+\#\{i:r_i<0\} \leq \dim_{\mathbb C}X-\#\{i:r_i\neq 0\}+\#\{i:r_i<0\}\leq \dim_{\mathbb C}X-\#\{i:r_i>0\}.
\]
For the virtual dimension constraint (\ref{vir-dim-3}) to hold, we must have
\[
\#\{i:r_i>0\}=0 \text{ and } [\gamma]_{-\vec r}=[\on{pt}]_{-\vec r}.
\]
Now the RHS of (\ref{iden-wdvv-mid-age}) is what we want.

For the LHS of (\ref{iden-wdvv-mid-age}), we first note that the markings with insertions $[\gamma]_{-\vec r}$ and $[\gamma^\vee]_{\vec r}$ must have mid-ages along the same divisors as mid-ages $\vec {\mathbf b}$.

Moreover, we must have $r_i=0$ for all other $i\in \{1,\ldots,n\}$ because $\cap_{i\in S_{\vec s}} D_{i}$ is a lowest dimensional stratum of $D$. We must have $\vec r=\vec {\mathbf b}+\vec l$ for some $\vec l\in \mathbb Z^n$. The virtual dimension constraint is stated in (\ref{vir-dim}). It can stated for invariants with mid-ages as follows: 

\begin{align*}
    \dim_{\mathbb C}X-3+3+0+0=0+0+d/2+\#\{i:r_i<0 \text{ or $r_i$ is mid-age}\},
\end{align*}
where $\gamma_j\in H^*(D_{\vec r})$ is of real degree $d$. As $D_{\vec{\mathbf b}}$ is a lowest dimensional stratum, we should have either $r_i=0$ or $r_i$ is a mid-age:
\[
\#\{i:r_i<0 \text{ or $r_i$ is mid-age}\}=\#\{i: \text{ $r_i$ is mid-age}\}.
\]
Therefore, the virtual dimensional constraint is
\begin{align*}
    \dim_{\mathbb C}X=d/2+(\dim_{\mathbb C}X-\dim_{\mathbb C}(D_{\vec{\mathbf b}})),
\end{align*}
where $\gamma_j\in H^*(D_{\vec{\mathbf b}})=H^*(D_{\vec{\mathbf b}+\vec l})$ is of real degree $d$. We need to have
\[
d/2=\dim_{\mathbb C}(D_{\vec{\mathbf b}}).
\]
Therefore, 
We must have $[\gamma]_{-\vec{\mathbf b}-\vec l}=[\on{pt}]_{-\vec{\mathbf b}-\vec l}$, for some $\vec l\in \mathbb Z^n$.

By Theorem \ref{thm-mid-age-mul}, the value of the invariants on the LHS: $\langle [1]_{\vec p},[1]_{\vec {\mathbf b}+\vec s},[\on{pt}]_{-\vec {\mathbf b}-\vec l}\rangle$ and $\langle [1]_{\vec {\mathbf b}+\vec l}, [1]_{\vec q}, [\on{pt}]_{-\vec {\mathbf b}}\rangle$ do not depend on the value of the mid-age. Set $\vec k=\vec s-\vec l$, we have
\[
\langle [1]_{\vec p},[1]_{\vec {\mathbf b}+\vec s},[\on{pt}]_{-\vec {\mathbf b}-\vec l}\rangle=\langle [1]_{\vec p},[1]_{\vec {\mathbf b}+\vec k},[\on{pt}]_{-\vec {\mathbf b}}\rangle
\]
where $\vec k+\vec l=\vec s\in \mathbb Z^n$. This concludes the proof.
\end{proof}

\if{
When $-\vec s\in B(\mathbb Z)$, we would like to show that the RHS is zero unless $\vec r=-\vec s$. This is because for mid-age invariants (for smooth pair first)
\[
\langle [1]_p,[1]_{\mathbf b},[\on{pt}]_{-r+\mathbf b}\rangle=0
\]
for $r>0$, unless $p=r$. Proof: apply the degeneration formula, all three markings have to map to the rubber. And there can not be edges, because the rubber invariant is zero unless degree zero and three markings, by the string equation. This is true when there are several negative contact orders.

\begin{corollary}
    The structure constant can be written in terms of mid-age invariants as follows
    \[
    N_{\vec p,\vec q,-\vec r}=\sum_{-\vec k-\vec l=\vec r\in B(\mathbb Z)}N_{\vec p, \vec{\mathbf b}+\vec k, -\vec{\mathbf{b}}}N_{\vec q, \vec{\mathbf b}+\vec l, -\vec{\mathbf{b}}}
    \]
    where the sum is over......
\end{corollary}
}\fi

\section{A comparison with logarithmic theta functions}\label{sec:comparison}

In this section, we use the result of \cite{BNR22}, \cite{BNR24} and \cite{Johnston} to explore the relation between our set-up and the set-up in intrinsic mirror symmetry \cite{GS19} where the structure constants are punctured logarithmic Gromov--Witten invariants. The notion of slope sensitivity was introduced in \cite{BNR22} to determine when log and orbifold invariants coincide.

\if{We consider the set-up in \cite{GS19} and \cite{GS21} but without assuming the divisor $D$ has a zero dimensional stratum. 

Let $P\subset H_2(X)$ be a finitely generated submonoid, containing all effective curve classes and the group of invertible elements $P^\times$ of $P$ coincides with the torsion part of $H_2(X)$. Let $I\subset P$ be a monoid ideal such that $P\setminus I$ is finite. We consider the multiplication rule
\begin{align}\label{multi-rule-finite}
\vartheta_{\vec p}\star \vartheta_{\vec q}=\sum_{\vec r\in B(\mathbb Z)} \sum_{\beta\in P\setminus I}N_{\vec p, \vec q, -\vec r}^{\beta}t^\beta\vartheta_{\vec r}.
\end{align}

Then there are only finitely many non-zero structure constants $N_{\vec p, \vec q, -\vec r}^{\beta}$ that we need to consider (\cite{GS19}*{Lemma 4.1}) and the sum defining $\vartheta_{\vec p}(x)$ is also finite (\cite{GS21}*{Lemma 3.22}). 
}\fi

The structure constants in (\ref{multi-rule-m}) can either be orbifold invariants or punctured log invariants. We use $N_{\vec p, \vec q, -\vec r}^{\on{orb},\beta}$ and $N_{\vec p, \vec q, -\vec r}^{\on{log},\beta}$ to denote orbifold structure constants and log structure constants respectively. Similarly, we use $\vartheta^{\on{orb}}_{\vec p}(\textrm x)$ and $\vartheta^{\on{log}}_{\vec p}(\textrm x)$ to denote orbifold theta functions and log theta functions.

Given $\vec p, \vec q, \vec r\in B(\mathbb Z)$ and an effective curve class $\beta\in H_2(X)$, we can consider a slope sensitive modification $\pi: \tilde X\rightarrow X$ with respect to this discrete data such that any lift of the discrete data to $\tilde X$ is slope sensitive. 

 Recall that the invariants defining $\vartheta_{\vec p}^{\on{orb}}$ are also three-point invariants and the pair of mid-age markings can be considered as a pair of markings with one sufficiently large contact order and one with sufficiently large negative contact order, where ``sufficiently large'' is with respect to the curve class $\beta$. 
 
 There are finitely many invariants in the WDVV equation (\ref{iden-wdvv-mid-age}), we may therefore pass to a sufficiently fine subdivision of $(X,D)$, $(\tilde X,\tilde D)$ such that each of the lift of the discrete data is slope sensitive. After potential further blow-up, we may assume that $\vec r$ is contained in a $1$-dimensional cone of $\Sigma(\tilde X)$. Then by \cite{Johnston}*{Corollary 7.4}, these orbifold structure constants of $(\tilde X,\tilde D)$ of degree $\tilde \beta$, agree with the corresponding log structure constants
 \[
 N_{\vec p, \vec q, -\vec r}^{\on{orb},\tilde \beta}=N_{\vec p, \vec q, -\vec r}^{\on{log},\tilde \beta}.
 \]
So we have structure constants in terms of three-point punctured invariants and the mid-age invariants in (\ref{iden-wdvv-mid-age}) are also three-point punctured invariants.

Up to this point, we did not say anything about the relation between the orbifold theta functions and the log theta functions.  Our theta functions are defined in terms of three-point orbifold invariants. When the discrete data is slope sensitive, the mid-age invariants agree with three-point punctured log invariants. On the other hand, invariants that define log theta functions in \cite{GS21} are in terms of two-point punctured invariants. 

Following \cite{Johnston}, for structure constants with contact order $-\vec r$ at the third marking, we can assume that, after blowing ups, $\vec r$ is in a $1$-dimensional cone of $\Sigma(\tilde X)$. That means, without loss of generality, we can assume the contact order is of the form $(-r,0,0,\ldots,0)$, where $r>0$. For the invariants in the theta functions, We also consider a further subdivision such that $\vec {\mathbf b}$ is in a $1$-dimensional cone of $\Sigma(\tilde X)$ and the contact orders of the second and the third markings are of the form $(b+k_1,0,\ldots,,0)$ and $(-b,-k_2,\ldots,-k_n)$, where $k_i\geq 0$ for $i\geq 2$. 


\begin{proposition}
    The following identity holds for three-point orbifold invariants
    \[
    \langle [1]_{\vec p},[1]_{(b+k_1,\vec 0)},[\on{pt}]_{(-b,-k_2,\ldots,-k_n)}\rangle_{0,3,\beta}^{(X,D)}=k_1\langle [1]_{\vec p}, [\on{pt}]_{(k_1,-k_2,-k_3,\ldots,-k_n)}\rangle_{0,2,\beta}^{(X,D)},
    \]
    where $k_i\geq 0$ for $i=2,\ldots,n$ and we use the convention that the RHS is zero if $k_1\leq 0$ (because the virtual dimension constraint does not satisfy for the invariant on the RHS if $k_1\leq 0$.)
\end{proposition}
\begin{proof}
The proof is similar to the proof of Proposition \ref{prop-mid-age-c}. We apply the degeneration to the normal cone with respect to $D_1$ and apply the degeneration formula to the invariant
\[
\langle [1]_{\vec p},[1]_{(b+k_1,\vec 0)},[\on{pt}]_{(-b,-k_2,\ldots,-k_n)}\rangle_{0,3,\beta}^{(X,D)}.
\]
The second marking $p_2$ and the third marking $p_3$ are distributed to the rubber integral over $D_1$ (with additional root structures). The rubber moduli can be pushed forward to the moduli space of stable maps to the orbifold-$D_1$. The insertion of the second marking $p_2$ becomes the identity class, so the rubber invariant is zero unless it is of degree zero and three markings. Therefore, the degeneration formula must give the invariant
\[
k_1\langle [1]_{\vec p}, [\on{pt}]_{(k_1,-k_2,-k_3,\ldots,-k_n)}\rangle\times 1,
\]
where multiplicity $k_1$ is coming from the edge of the degeneration graph and the rubber invariant is $1$, the insertion of the second marking is of the form $[\on{pt}]_{(k_1,-k_2,-k_3,\ldots,-k_n)}$ because the rubber moduli has to be of degree zero and three markings. 
\end{proof}

\begin{remark}
    The canonical wall structure defined in \cite{GS21} uses punctured invariants with the assumption that the divisor $D$ has a zero dimensional stratum. Our results indicate that there should also be the canonical wall structure using punctured or orbifold invariants without assuming $D$ has a zero dimensional stratum. The fact that there are well-defined theta functions indicate that the canonical wall structure is consistent. We believe that this generalization is expected by experts. However, without the assumption of a zero dimensional stratum, this structure does not give a mirror construction. One motivation to study this generalization is the mirror duality between degenerations and fibrations as we explained in the introduction.
\end{remark}

\bibliographystyle{amsxport}
\bibliography{main}

\end{document}